\documentclass[11pt, a4paper]{article}
\usepackage[latin1]{inputenc}
\usepackage[english]{babel}
\usepackage{indentfirst}
\usepackage{amstext}
\usepackage{amsfonts}
\usepackage{textcomp}
\usepackage{amssymb}
\usepackage{amscd}
\usepackage[dvips]{graphicx}
\usepackage{setspace}
\usepackage{fancyhdr}

\newcommand {\demo}{\hskip -0.6cm{\bf Proof.  }}
\newcommand {\fim}{\hfill{$\square$}\vskip 1pc}
\newcommand {\nl}{\newline}
\newcommand {\cl}{\centerline}
\newcommand {\R}{\mathbb{R}}
\newcommand {\N}{\mathbb{N}}
\newcommand {\Z}{\mathbb{Z}}

\newcommand {\F}{\mathbb{F}}

\newcommand {\la}{\lambda}
\newcommand {\al}{\alpha}
\newcommand {\til}{\widetilde}
\newcommand {\lb}{\linebreak}

\newcommand {\ovl}{\overline}
\newcommand {\m}[1]{{\widetilde{#1}\,}}
\newcommand{\af}[2]{\vspace{0.3cm}\hskip -0.6cm {\it Claim#1: #2}\vspace{0.3cm}}
\newcommand{\me}[1]{(#1)^{\hspace{-0.1cm}\widetilde{\hspace{0.3cm}}}}
\newcommand{\supp}{\text{supp}}
\newcommand{\tnorm}{\|\hspace{-0.03cm}|}
\newcommand{\fimaf}{\vspace{0.3cm}}

\newcommand {\prcr}[1]{\mathcal{O}(#1,\al,L)}
\newcommand {\tp}[1]{\mathcal{T}(#1,\al,L)}
\newcommand{\funcao}[5]{\begin{array}{lrcl}
#1:&\!\!\!#2 & \rightarrow & #3 \\
  &\!\!\! #4 & \mapsto & #5
\end{array}}

\newtheorem{teorema}{Theorem}[section]
\newtheorem{lema}[teorema]{Lemma}
\newtheorem{corolario}[teorema]{Corollary}
\newtheorem{definicao}[teorema]{Definition}
\newtheorem{proposicao}[teorema]{Proposition}
\begin{document}
\onehalfspace

\cl{\large THE CROSSED PRODUCT BY A PARTIAL ENDOMORPHISM}
\vspace{1cm}
\cl{R. Exel\footnote{Partially supported by Cnpq} and D. Royer\footnote{Supported by Cnpq}}
\vspace{1cm}
\begin{abstract}

Given a closed ideal $I$ in a $C^*$-algebra $A$, an ideal $J$ (not necessarily closed) in $I$, a
*-homomorphism $\al:A\rightarrow M(I)$ and a map 
$L:J\rightarrow A$ with some properties, based on
$\cite{pimsner}$ and
$\cite{katsura}$ we define a $C^*$-algebra $\prcr{A}$ which we call the {\it Crossed Product by a Partial
Endomorphism.}

In the second section we introduce the Crossed Product by a Partial Endomorphism $\prcr{X}$ induced by a
local homeomorphism $\sigma:U\rightarrow X$ where
$X$ is a compact Hausdorff space and $U$ is an open subset of $X$.
The
main result of this section is that every nonzero gauge invariant ideal of $\prcr{X}$ has nonzero intersection
with $C(X)$. We present the example which
motivated this work, the Cuntz-Krieger algebra for infinite matrices (see \cite{exelmatinf}).

We show in the third section a bijection between the gauge invariant ideals of $\prcr{X}$ and the
$\sigma,\sigma^{-1}$-invariant open subsets of $X$.

The last section is dedicated to the study of $\prcr{X}$ in the case where the pair $(X,\sigma)$ 
has an extra property, wich we call topological freeness. We prove that in this case every nonzero
ideal of $\prcr{X}$ has nonzero intersection with $C(X)$. If moreover $(X,\sigma)$ has the property that
$(X',\sigma_{|_{X'}})$ is topologically free for each closed $\sigma,\sigma^{-1}$-invariant subset $X'$ of $X$
then we obtain a bijection between the ideals of $\prcr{X}$ and the open
$\sigma,\sigma^{-1}$-invariant subsets of $X$. We conclude this section by showing a simplicity criteria for
the Cuntz-Krieger algebras for infinite matrices.
\end{abstract}

\tableofcontents

\section*{Introduction}\addcontentsline{toc}{section}{Introduction}

In \cite{exelprodcruz} it was introduced by the first named author the concept of Crossed Product by an Endomorphism,
based on a C$^*$-dynamical system $(A,\al,L)$.
In this article it was shown that the Cuntz-Krieger algebra is an example of Crossed Product by an Endomorphism. 
The $C^*$-dynamical system associated to this example is induced by the Markov subshift $(\Omega_A,\sigma)$,
that is, the endomorphism $\al:C(\Omega_A)\rightarrow C(\Omega_A)$ is given by $\al(f)=f\circ \sigma$ and 
$L:C(\Omega_A)\rightarrow C(\Omega_A)$ is defined by \lb
$L(f)(x)=\frac{1}{\#\sigma^{-1}(x)}\sum\limits_{y\in \sigma^{-1}(x)}f(y)$
for each $x\in X$ and for each $f\in C(\Omega_A)$.

It was defined in \cite{exelmatinf} by the first named author and M. Laca the 
Cuntz-Krieger algebra for infinite matrices. This algebra has a topological compact Hausdorff space $\m{\Omega_A}$
associated to it, which can be seen in [\ref{exelmatinf}:4-7]. The difference between this case and the
previous one is that the shift
$\sigma$ can not be defined in the whole space $\m{\Omega_A}$, but only in an open subset $U$ of $\m{\Omega_A}$. 
Then the local homeomorphism $\sigma:U\rightarrow \m{\Omega_A}$ induces the *-homomorphism
$\al:C(\m{\Omega_A})\rightarrow C^b(U)$ given by $\al(f)=f\circ\sigma$. Moreover, since $\#\sigma^{-1}(x)$ may
be infinite for some $x\in \m{\Omega_A}$, the convergence of the sum $\sum\limits_{y\in
\sigma^{-1}}f(y)$ is not guaranteed and so $L(f)$ can not be defined by 
$L(f)(x)=\sum\limits_{y\in \sigma^{-1}}f(y)$ for every $f\in C(\m{\Omega_A})$.
However, we will show that for each $f\in C_c(U)$, $L(f)$ defined by
$L(f)(x)=\sum\limits_{y\in \sigma^{-1}(x)}f(y)$ for each $x\in \m{\Omega_A}$ is an element of $C(\m{\Omega_A})$. 
In this way we obtain a map
$L:C_c(U)\rightarrow C(\m{\Omega_A})$. Because $\al$ is not an endomorphism in $C(\m{\Omega_A})$ 
and the domain of $L$ is not the whole algebra $C(\m{\Omega_A})$, the triple
$(A,\al,L)$ (which we also call by $C^*$-dynamical system) is not a $C^*$-dynamical system as in \cite{exelprodcruz} 
and therefore the construction of Crossed Product by an Endomorphism defined in \cite{exelprodcruz} cannot be applied.

In this work we define, making use of the constructions of T. Katsura (\cite{katsura})
and M. Pimsner (\cite{pimsner}), the {\it Croseed Product by a Partial Endomorphism.} We show that  
our construction may be applied to the situation described in the previous paragraph. We study specially the
case where the Crossed Product by a Partial Endomorphism, wich we denote by $\prcr{X}$, is induced
by a local homeomorphism $\sigma:U\rightarrow X$, where $U$ is an open subset of a topological compact
Hausdorff space $X$. More specifically, we show a bijection between the gauge invariant ideals of $\prcr{X}$
and the $\sigma,\sigma^{-1}$-invariant open subsets of $X$. Moreover, if $(X,\sigma)$ has the property that
$(X',\sigma_{|_{X'}})$ is topologically free for
every closed $\sigma,\sigma^{-1}$-invariant subset $X'$ of $X$ then there exists a bijection between the
ideals of $\prcr{X}$ and the open $\sigma,\sigma^{-1}$-invariant subsets of $X$. Finally we present a
simplicity criteria for the Cuntz-Krieger algebras for infinite matrices.

The choice of the name Crossed Product by a Partial Endomorphism for the algebra $\prcr{A}$ defined in this work
was motived by the local homeomorphism $\sigma:U\rightarrow X$ where $U$ is an open subset of $X$.

In \cite{kwasniewski}, B. K. Kwasniewski defined an algebra which he called {\it Covariance algebra of a
partial dynamical system} based on a partial dynamical system $(X,\al)$, that is, a continuous map
$\al:\Delta\rightarrow X$ where $X$ is a compact Hausdorff space and $\Delta$ is a clopen subset of $X$ and
$\al(\Delta)$ is open. In our construction $\Delta$ need not be clopen, only open, but we require that $\al$
is a local homeomorphism. The possible relationship between these two constructions will be studied in a
future paper.  

\section{The Crossed Product by a Partial Endomorphism}

In this section we define the crossed product by a partial
endomorphism and show some results about its structure. We study the gauge action and gauge-invariant ideals
of this algebra.

\subsection{Definitions and basic results} Let $A$ be a $C^*$-algebra and $I$ a closed two-sided ideal in $A$.

\begin{definicao}
A partial endomorphism is a *-homomorphism $\al:A\rightarrow M(I)$ where $M(I)$ is the multiplier algebra of
$I$.
\end{definicao}

Let $J$ be a two-sided self adjoint idempotent (not necessarily closed) ideal in
$I$ and let 
$\al:A\rightarrow M(I)$ and $L:J\rightarrow A$ be functions. We denote a such situation by $(A, \al, L)$.

\begin{definicao}
$(A,\al,L)$ is a $C^*$-dynamical system if $(A,\al,L)$ has the following properties:
\begin{itemize}

\item $\al$ is a partial endomorphism, \item $L$ is linear, positive and
preserves *, \item $L(\al(a)x)=aL(x)$ for all $a$ in $A$ and $x$
in $J$.
\end{itemize}
\end{definicao}

The function $L$ is positive in the sense that $L(x^*x)$ is a
positive element of $A$ for all $x$ in $J$. Moreover, denoting $\al(a)$ by $(L^a, R^a)$, $\al(a)x$ is a
notation for the element $L^a(x)$. Note that if $x,y\in J$ and $a\in A$ then $L^a(x)\in I$ and so
$L^a(xy)=L^a(x)y\in J$. Since $J$ is idempotent we have in general that
$\al(a)x\in J$ for all $a\in A$ and $x\in J$. Therefore $\al(a)x$ lies in fact in the
domain of $L$. Defining $x\al(a)=R^a(x)$ for all $x\in J$ and $a\in A$ we have that $(\al(a)x)^*=x^*\al(a^*)$ 
for every $x\in J$ and $a\in A$. In fact,
$(\al(a)x)^*=(L^a(x))^*=(R^a)^*(x^*)=R^{a^*}(x^*)=x^*\al(a^*).$
In the same way $(x\al(a))^*=\al(a^*)x^*$. 

If $(A,\al,L)$ is a $C^*$-dynamical system then $L(x\al(a))=L(x)a$ for all $a\in A$ and
$x\in J$. In fact, given $a\in A$ e $x\in J$, since $a^*\in A$ and $x^*\in J$ we have that
$L(\al(a^*)x^*)=a^*L(x^*)$. Therefore
$L(x\al(a))=L((x\al(a))^*)^*=L(\al(a^*)x^*)^*=(a^*L(x^*))^*=L(x)a.$

The next goal is to define a left $A$-module which is also a right Hilbert $A$-module.
Define the operation
$$\funcao{.}{J\times A}{J}{(x,a)}{x\al(a)}.$$ It is easy to verify that this operation is
bilinear and associative. Thus $J$ is a right $A$-módule. It is also easy to see that the
function
$$\funcao{\langle ,\rangle }{J\times J}{A}{(x,y)}{L(x^*y)}$$
is a semi-inner product. Considering the quotient of $J$ by $N_0=\{x\in
J:\langle x,x\rangle =0\}$ and denoting the elements $x$ of $J$ by $\m{x}$ in $J/N_0$ (or
by $\me{x}$) we obtain an inner product of $J/N_0$ in $A$ defined by $\langle
\m{x},\m{y}\rangle =\langle x,y\rangle $. So the function
$$\funcao{\| \,\|}{J/N_0}{\R^+}{\m{x}}{\sqrt{\|\langle
\m{x},\m{x}\rangle \|}}$$ defines a norm in $J/N_0$. Denote by $M$ the right Hilbert
$A$-module $\ovl{(J/N_0)}^{\| \|}$.

Let us now define a left $A$-module structure for $M$. Given  $a\in A$ and $x\in J$ we
have that $x^*a^*ax,\|a\|^2x^*x\in J$. Since $x^*(\|a\|^2-a^*a)x$ may be written in the
form $(bx)^*(bx)$ with $bx\in J$ we have that $L(x^*\|a\|^2x-x^*a^*ax)\geq 0$ and so
$L(x^*a^*ax)\leq \|a\|^2L(x^*x)$ from where $\|L(x^*a^*ax)\|\leq \|a\|^{2}\|L(x^*x)\|$.
Therefore
$$\|\m{ax}\|^2=\|\langle \m{ax},\m{ax}\rangle \|=\|L(x^*a^*ax)\|\leq
\|a\|^2\|L(x^*x)\|=\|a\|^2\|\langle \m{x},\m{x}\rangle \|=\|a\|^2\|\m{x}\|^2,$$ and so,
$\|\m{ax}\|\leq \|a\|\|\m{x}\|$. This allows us define the operation
$$\funcao{.}{A\times M}{M}{(a,m)}{am},$$ where $a\m{x}=\m{ax}$, which is
bilinear and associative, and so $M$ is a left $A$-module. This operation gives rise to a
*-homomorphism from $A$ in $L(M)$. In fact, defining $\varphi:A\rightarrow L(M)$ by
$\varphi(a)m=am$ we have:

\begin{proposicao}
$\varphi$ is a *-homomorphism.
\end{proposicao}

\demo For all $a\in A$, $\varphi(a):M\rightarrow M$ definded by $\varphi(a)(m)=am$ for all $m\in M$ is a linear function.
Moreover, for $x,y\in J$,
$$\langle \varphi(a)\m{x},\m{y}\rangle =\langle
\m{ax},\m{y}\rangle =L((ax)^*y)=L(x^*a^*y)=\langle \m{x},\m{a^*y}\rangle =\langle
\m{x},\varphi(a^*)\m{y}\rangle,$$ and since $J/{N_0}$ is dense in $M$ it follows that
$\langle \varphi(a)m,n\rangle =\langle m,\varphi(a^*)n\rangle $ for all $m,n\in M$. This
shows that $\varphi(a)$ is adjointable and $\varphi(a)^*=\varphi(a^*)$. Obviously
$\varphi$ is linear and multiplicative. \fim

\begin{definicao}
The Toeplitz algebra $\tp{A}$ associated to the $C^*$-dynamical system $(A,\al,L)$ is the
universal $C^*$-algebra generated by $A\cup M$ with the relations of $A$, of $M$, the
$A$-bi-module products and $m^*n=\langle m,n\rangle $ for all $m,n\in M$.
\end{definicao}

Note that the universal algebra in fact exists, since the relations are admissible. We
will denote by $\widehat{K_1}$ the closed sub-algebra of $\tp{A}$ generated by the
elements of the form $mn^*$, for $m,n\in M$.

\begin{definicao}
A redundancy in $\tp{A}$ is a pair $(a,k)$ where $a\in A$, $k\in \widehat{K_1}$ and
$am=km$ for all $m\in M$.
\end{definicao}

Let $I_0=\text{ker}(\varphi)^\bot\cap \varphi^{-1}(K(M))$
where $\varphi:A\rightarrow L(M)$ is the *-homomorphism given by the left
multiplication.

\begin{definicao}
The Crossed Product by a partial Endomorphism associated to the $C^*$-dynamical system $(A,\al,L)$ is the
quotient of $\tp{A}$ by the ideal generated by the elements $a-k$ for all redundancies
$(a,k)$ such that $a\in I_0$, and will be denoted by $\prcr{A}$.
\end{definicao}

It follows from \cite{katsura} that $A\ni a\rightarrow a\in \prcr{A}$ is injective. In
the following proposition will be showed some consequences of this fact. Let us temporarily denote by $\widehat{a}$ and $\widehat{m}$ the elements of $A$ and $M$ in $\tp{A}$.
Define
$$\widehat{K_n}=\ovl{\text{span}}\{\widehat{m_1}\cdots\widehat{m_n}\widehat{l_1}^*\cdots\widehat{l_n}^*:\,\,m_i,l_i\in
M\}$$ and denote by $q$ the quotient map from $\tp{A}$ to $\prcr{A}$.

\begin{proposicao}\label{consequencias}
\begin{itemize}
\item[a)] $A\ni a\mapsto q(\widehat{a})\in \prcr{A}$ is an injective *-homomorphism.

\item[b)] $A\ni a\mapsto \widehat{a}\in \tp{A}$ and $q_{|_{\widehat{A}}}$ are injective
*-homomorphisms.

\item[c)] $M\ni m\mapsto \widehat{m}\in \tp{A}$ is an isometry.

\item[d)] $q_{|_{\widehat{M}}}$ is an isometry.

\item[e)] $M\ni m\mapsto q(\widehat{m})\in \prcr{A}$ is an isometry.

\item[f)] $q_{|_{\widehat{K_n}}}$ is an injective *-homomorphism.

\end{itemize}
\end{proposicao}

\demo a) Is a consequence of \cite{katsura}\nl b) Follows from a)\nl c) Given $m\in M$,
$\|\widehat{m}\|^2=\|\widehat{m}^*\widehat{m}\|=\|\widehat{\langle
m,m\rangle }\|.$ Since $\langle m,m\rangle \in A$, it follows from b) that
$\|\widehat{\langle m,m\rangle }\|=\|\langle m,m\rangle \|$. Moreover
$\|m\|^2=\|\langle m,m\rangle \|$. Then $\|\widehat{m}\|^2=\|\langle
m,m\rangle \|=\|m\|^2.$ \nl d) For all $\widehat{m}\in \widehat{M}$ we have
$\widehat{m}^*\widehat{m}\in \widehat{A}$. By a), $q_{|_{\widehat{A}}}$ is injective and
therefore an isometry. Then
$\|q(\widehat{m})\|^2=\|q(\widehat{m}^*\widehat{m})\|=\|\widehat{m}^*\widehat{m}\|=\|\widehat{m}\|^2.$\nl
e) Follows from c) and d).\nl f) Let $k\in \widehat{K_n}$ and suppose $q(k)=0$. Then
$q((\widehat{M}^*)^nk\widehat{M}^n)=0$. Since $(\widehat{M}^*)^nk\widehat{M}^n\subseteq
\widehat{A}$ it follows from b) that $(\widehat{M}^*)^nk\widehat{M}^n=0$. Then
$\widehat{K_n}k\widehat{K_n}=0$ and so $k=0$. \fim

From now on we will identify the elements $\widehat{a}\in\tp{A}$ and $q(\widehat{a})\in \prcr{A}$ with the element $a$
of $A$. This notation will not cause confusion, by a) and b) of the previous
proposition. In the same way, justified by c) and e) we will identify the elements $\widehat{m}\in \tp{A}$ and
$q(\widehat{m})\in\prcr{A}$ with the element $m\in M$. With these identifications,
$$\widehat{K_n}=\ovl{\text{span}}\{m_1\cdots m_nl_1^*\cdots l_n^*:\,\,m_i,l_i\in
M\}\subseteq \tp{A}.$$ Define
$$K_n=\ovl{\text{span}}\{m_1\cdots m_nl_1^*\cdots l_n^*:\,\,m_i,l_i\in
M\}\subseteq \prcr{A}$$ and note that $q(\widehat{K_n})=K_n$. If $(a,k)\in A\times
\widehat{K_1}$ is a redundancy and $a\in I_0$ then $q(a)=q(k)$. Since $a=q(a)$ in
$\prcr{A}$ it follows that $a=q(k)$ in $\prcr{A}$.

The spaces $K_n$ e $\widehat{K_n}$ are clearly closed under the sum and are self-adjoint.
Moreover, the following proposition shows that they are closed under multiplication, and
so are $C^*$-algebras.

\begin{proposicao}\label{subalgebras}
a) $\widehat{K_n}\widehat{K_m}\subseteq \widehat{K_{\text{max}\{n,m\}}}$ and also
$K_nK_m\subseteq K_{\text{max}\{n,m\}}$. \nl b) $A\widehat{K_n}\subseteq\widehat{K_n}$,
$\widehat{K_n}A\subseteq\widehat{K_n}$ and also $AK_n\subseteq K_n$ and $K_nA\subseteq
K_n$.
\end{proposicao}

\demo Since $K_n=q(\widehat{K_n})$ it suffices to show the result for the algebra
$\tp{A}$.\nl a) Taking adjoins we may suppose $n\leq m$. Given $l_1...l_nt_1^*...t_n^*\in
\widehat{K_n}$ and $p_1...p_mq_1^*...q_m^*\in \widehat{K_m}$, how
$a=t_1^*...t_n^*p_1...p_n\in A$ it follows that $l_na\in M$. Therefore
$$l_1...l_nt_1^*...t_n^*p_1...p_mq_1^*...q_m^*=l_1...l_nap_{n+1}...p_mq_1^*...q_
m^*\in \widehat{K_m}.$$ This is enough since $\widehat{K_n}$ are generated by elements of
this form.\nl b) Follows by the fact that $am\in M$ for all $a\in A$ and $m\in M$.\fim

We will denote by $m\otimes n$ the element of $K(M)$ given by $m\otimes n(\xi)=m\langle
n,\xi\rangle$, for all $\xi\in M$.

\begin{proposicao}\label{pe1} There
exists a *-isomorphism
$S:\widehat{K_1}\rightarrow K(M)$ such that $S(mn^*)=m\otimes n.$
\end{proposicao}

\demo Given $k\in \widehat{K_1}$ and $m\in M$ then $km\in M$ because $M$ is closed in
$\tp{A}$ by the proposition \ref{consequencias} c). In $\tp{A}$,
$\langle km,n\rangle =(km)^*n=m^*k^*n=\langle m,k^*n\rangle,$ and how
$\langle m,k^*n\rangle,\langle km,n\rangle\in A$, by \ref{consequencias} b) $\langle
m,k^*n\rangle=\langle km,n\rangle$ in $A$. So, defining
$S(k):M\rightarrow M$ by $S(k)(m)=km$ it follows that
$\langle S(k)m,n\rangle=\langle km,n\rangle=\langle m,k^*n\rangle=\langle
m,S(k^*)n\rangle$ for all $m,n\in M$. This shows that $S(k)$ is adjointable and
$S(k)^*=S(k^*)$. Since $S(k)\in L(M)$ we may define $S:\widehat{K_1}\rightarrow L(M)$ which is clearly
linear and multiplicative, and so $S$ is a *-homomorphism. Obviously $S(mn^*)=m\otimes
n$, and therefore $S(k)\in K(M)$ for all $k\in \widehat{K_1}$. Moreover
$S(\widehat{K_1})$ is a dense set in $K(M)$ and so $S(\widehat{K_1})=K(M)$. In order to see that
$S$ is injective suppose $S(k)=0$, that is, $kM=0$. Then $k\widehat{K_1}=0$ and since $k\in
\widehat{K_1}$ it follows that $k=0$.\fim

If $(a,k)$ is a redundancy then $am=km$ for all $m\in M$, from where
$\varphi(a)(m)=S(k)(m)$ for each $m\in M$. Since $S(k)\in K(M)$ it follows that $\varphi^{-1}(a)\in K(M)$.
So the algebra $\prcr{A}$ coincides with the quotient of $\tp{A}$ by the
ideal generated by the elements of the form $(a-k)$ for all redundancy $(a,k)$ such that
$a\in \text{ker}(\varphi)^\bot$.

Given a $C^*$-dynamical system $(A,\al,L)$ and a closed ideal $N$ in $A$ such that
$J\subseteq N\subseteq I$, we may consider an other $C^*$-dynamical system $(A,\beta,L)$
where the partial endomorphism \lb $\beta:A\rightarrow M(N)$ is given by
$\beta(a)=(L_{|_N}^a,R_{|_N}^a)$, considering that $\al(a)=(L^a,R^a)$. Since\lb
$x\beta(a)=x\al(a)$ for all $x\in J$ and $a\in A$ it follows that
$\prcr{A}=\mathcal{O}(A,\beta,L)$. By this reason we may suppose that $J$ is a dense
ideal in $I$. This situation will occur in the second section.

It may be showed without much difficulty that the crossed product by endomorphism
introduced in \cite{exelprodcruz} in some situations may be seen as crossed products by
a partial endomorphism. More specifically, this holds if $\langle \al(A)\rangle=A$ and
$L$ is faithfull or if $\al:A\rightarrow A$ is injective, $\al(A)=\al(1)A\al(1)$, and
$L:A\rightarrow A$ is given by $L(a)=\al^{-1}(\al(1)a\al(1))$. The first situation occurs in 
Cuntz-Krieger algebras (see [\ref{exelprodcruz}: 6]) end the last situation occurs in 
Pashke's crossed product and in the crossed product proposed by Cuntz (see [\ref{exelprodcruz}]).

\subsection{The gauge action}\label{ck}

The next goal is to show that every gauge-invariant ideal of $\prcr{A}$ has non-trivial intersection with the fixed point algebra of the gauge action in $\prcr{A}$.

By the universal property of $\tp{A}$ it follows that for each $\la \in S^1$ there exists
a \lb*-homomorphism
$\theta_\la:\tp{A}\rightarrow\tp{A}$ which satisfies $\theta_\la(a)=a$ for all $a$ in $A$ and $\theta_\la(m)=\la m$
for all $m\in M$.
If $(a,k)$ is a redundancy, because $\theta_\la(a)=a$ and $\theta_\la(k)=k$ it follows that
$(\theta_\la(a),\theta_\la(k))$ is also a redundancy, and so we may consider
$\theta_\la:\prcr{A}\rightarrow \prcr{A}$. Note that
$\theta_{\la_1}\theta_{\la_2}=\theta_{\la_1\la_2}$ from where $\theta_\la$ is a
*-automorphism, with inverse $\theta_{\overline{\la}}$. Moreover, given $r\in \prcr{A}$,
the function
$S^1\ni \la\mapsto \theta_\la(r)\in \prcr{A}$ is continuous. Then we may consider
$$\funcao{E}{\prcr{A}}{\prcr{A}}{r}{\int\limits_{S^1}\theta_\la(r)d\la}.$$

\begin{proposicao}
The fixed point algebra of $\theta$ is $K=\overline{span}\{A,K_n;n\in \N\}$ and $E$ is a
faithful conditional expectation onto $K$.
\end{proposicao}

\demo It is not difficult to show that $E$ is a faithful conditional expectation onto the
fixed point algebra. So it suffices to show that $\text{Im}(E)=K$. The equality holds because
$$E(am_1\cdots m_kn_1^*\cdots n_l^*b)= \left\{\begin{array}{ll}
am_1\cdots m_kn_1^*\cdots n_l^*b & \text{ se } k=l\\ 0 & \text{ se
}k\neq l
\end{array}\right.$$ and the space generated by elements of the form $am_1\cdots m_in_1^*\cdots m_j^*b$ is dense in $\prcr{A}$.\fim

\begin{definicao}\label{defgaugeinv}
A ideal $I$ in $\prcr{A}$ is gauge-invariant if $\theta_\la(I)\subseteq I$ for each
$\la\in S_1$.
\end{definicao}

If $I$ is gauge-invariant, the gauge action in $\prcr{A}/I$ is given by
$$\funcao{\beta_\la}{\prcr{A}/I}{\prcr{A}/I}
{\pi(r)}{\pi(\theta_\la(r))},$$ where $\pi$ is the quotient map. In this case $\pi$ is
covariant by the gauge actions $\theta$ and $\beta$, in the sense that
$\pi(\theta_\la(r))=\beta_\la(\pi(r))$ for all $r\in \prcr{A}$ and for each $\la\in S^1$.
Moreover, the fixed point algebra for $\beta$ is $\pi(K)$ because the conditional
expectation $F$ induced by $\beta$ is such that $F(\pi(r))=\pi(E(r))$ for each $r\in \prcr{A}$.

\begin{proposicao}\label{idintk}
If $0\neq I\unlhd \prcr{A}$ is gauge-invariant then $I\cap K\neq 0$.
\end{proposicao}

\demo Since $\theta_\la(I)\subseteq I$ for all $\la\in S^1$ then $E(r)\in I$ for all $r\in
I$. By the fact that $E$ is faithful it follows that, given $0\neq r\in I$ then $E(r^*r)\neq 0$. Since $E(r^*r)\in
K\cap I$, the result is proved.
\fim

Defining $$L_0=A \text{ and } L_n=A+K_1+\cdots+K_n  \text{ for every } n\geq 1$$ we have that
$L_0\subseteq L_1\subseteq L_2\subseteq \cdots$ and
$K=\ovl{\bigcup\limits_{n\in \N}L_n}.$ This form to see the algebra $K$ will be useful in some situations
which will appear latter. In some of this situations we will use the fact, given by the
following proposition, that the algebras $L_n$ (by the proposition \ref{subalgebras}
$L_n$ are algebras) are closed, for all $n\in \N$.

\begin{proposicao}
For each $n\in \N$ the algebras $L_n$ are closed.
\end{proposicao}

\demo The case $L_0$ follows by \ref{consequencias} a). By induction suppose $L_n$
closed. Note that $K_{n+1}\unlhd \ovl{L_{n+1}}$ and that $L_n$ is a closed sub-algebra of
$\ovl{L_{n+1}}$. By [\ref{pedersen}: 1.5.8], $L_n+K_{n+1}$ is a closed sub-algebra of
$\ovl{L_{n+1}}$. Therefore $L_{n+1}=L_n+K_{n+1}=\ovl{L_n+K_{n+1}}=\ovl{L_{n+1}}$. \fim

\section{The Crossed Product by a Partial Endomorphism induced by a local homeomorphism}

Given a topological compact Hausdorff space $X$ and a local homeomorphism
$\sigma:X\rightarrow X$, defining
$\al:C(X)\rightarrow C(X)$ by $\al(f)=f\circ\sigma$ and $L:C(X)\rightarrow C(X)$
by $L(f)(x)=\sum\limits_{y\in \sigma^{-1}(x)}f(y)$ for all $x\in X$, we obtain a $C^*$-dynamical system.
This situation occurs in the Cuntz-Krieger algebra in \cite{exelprodcruz}. A more general
situation consists in considering an open set $U\subseteq X$ and a local homeomorphism
$\sigma:U\rightarrow X$. In this case, defining $\al$ as above, for all $f\in C(X)$
$\al(f)$ is an element of $C^b(U)$, where $C^b(U)$ is the set of all continuous and
bounded functions in $U$. Moreover, $\#\sigma^{-1}(x)$ may be infinite for some $x\in X$,
and therefore $L$ can not be defined as above.

Although, if $f\in C_c(U)$, that is, $f\in C(X)$ such that
$\text{supp}(f)=\ovl{\{x\in X:f(x)\neq 0\}}\subseteq U$, we will show that
$\sum\limits_{y\in \sigma^{-1}(x)}f(y)$ involves finitely many summands for every $x\in X$. We will also show
that, for each $f\in C_c(U)$, $L(f)$ defined by $L(f)(x)=\sum\limits_{y\in
\sigma^{-1}(x)}f(y)$ is an element in $C(X)$, and so we may define $L:C_c(U)\rightarrow
C(X)$. Moreover, since $C^b(U)$ and $M(C_0(U))$ are *-isomorphic we obtain a partial
endomorphism $\widetilde{\al}:C(X)\rightarrow M(C_0(U))$.

We begin this section by showing that $(C(X),\widetilde{\al},L)$ is a $C^*$-dynamical
system which will give us the crossed product by a partial endomorphism $\prcr{X}$.

The second part is dedicated to presenting some basic results about the structure of
$\prcr{X}$, and the most important result of this part is that every ideal of
$\prcr{X}$ which has nonzero intersection with $K$ (the fixed point
algebra of the gauge action) has nonzero intersection with $C(X)$.

In the last part we show that the Cuntz-Krieger algebra for infinite
matrices (see \cite{exelmatinf}) is a crossed product by a partial endomorphism.
This is the example which motivated this work.

The choice of the name {\it Crossed Product by a Partial Endomorphism} for the algebra $\prcr{A}$
was motivaded by the local homeomorphism $\sigma$.

\subsection{The algebra $\prcr{X}$}\label{dynamicalsystem}

Let $X$ be a topological compact Hausdorff space, $U\subseteq X$ an open subset and
$\sigma:U\rightarrow X$ a local homeomorphism.
Define $$\begin{array}{cc}\al:C(X)\rightarrow C^b(U)\\
f\mapsto f\circ\sigma
\end{array}$$
which is a *-homomorphism. For each $f\in C_c(U)$ define for all $x\in X$,
$$L(f)(x)=\left\{\begin{array}{ll} \sum\limits_{y\in U\atop\sigma(y)=x}f(y)&\text{ if }
\sigma^{-1}(x)\neq \emptyset\\ 0&\text{ otherwise}
\end{array}\right ..$$

If $K\subseteq U$ is a compact subset, taking an open cover $U_1,\cdots ,U_n$ of $K$ in
$U$ such that $\sigma_{|_{U_i}}$ is homeomorphism, for every $x\in X$ there exists no
more than one element $x_i$ in each $\sigma^{-1}(x)\cap U_i$. Therefore there exists at
most $n$ elements in $\sigma^{-1}(x)\cap K$. It follows that the sum which defines
$L(f)(x)$ involves finitely many summands for each $x\in X$, and so $L(f)(x)$ in fact may be defined as above.

\begin{lema}
For each $f\in C_c(U)$, $L(f)$ is an element of $C(X)$.
\end{lema}

\demo Let $f\in C_c(U)$ and $K=\supp(f)$. We will show that $L(f)$ is continuous on each point of $X$.
Given $x\in X\setminus \sigma(K)$, since $X\setminus \sigma(K)$ is open and
$L(f)y=0$ for all $y\in X\setminus\sigma(K)$, it follows that $L(f)$ is
continuous in $x$. Let $x\in \sigma(K)$, $\{x_1,\cdots ,x_k\}=\sigma^{-1}(x)\cap K$,
and $U_j$ open disjoint neighbourhoods of $x_j$ such that $\sigma_{|_{U_j}}$ is
a homeomorphism. The $U_j$ may be taken such that $\sigma(U_j)$ are open,
because $\sigma$ is a local homeomorphism.

\af{}{There exists an open set $V\ni x$ such that $\sigma^{-1}(V)\cap (K\backslash
(\bigcup\limits_{j=1}^kU_j))=\emptyset$.}

Suppose $\sigma^{-1}(V)\cap (K\backslash
(\bigcup\limits_{j=1}^kU_j))\neq\emptyset$ for each open set $V$ which contains
$x$. For every open subset $W\ni x$ define
$$F_W=\sigma^{-1}(\overline{W})\cap(K\backslash
(\bigcup\limits_{j=1}^kU_j)).$$ Since $\sigma^{-1}(\ovl{W})$ is closed in
$U$ and $K\backslash (\bigcup\limits_{j=1}^kU_j)\subseteq U$ is compact, it
follows that $F_W$ is compact, and therefore closed in $X$. Moreover $F_W$ is
nonempty because
$$\emptyset \neq \sigma^{-1}(W)\cap (K\backslash
(\bigcup\limits_{j=1}^kU_j))\subseteq F_W.$$ Given $W_1,..., W_m$ open
neighbourhoods of $x$, we have that $F_{\bigcap\limits_{j=1}^mW_j}\subseteq F_{W_j}$ for
each $j$ from where
$F_{\bigcap\limits_{j=1}^mW_j}\subseteq \bigcap_{j=1}^mF_{W_j},$ and so
$\bigcap\limits_{j=1}^m F_{W_j}\neq\emptyset$ for each finite collection of
open neighbourhoods $W_1,...,W_m$ of $x$. By the fact that $X$ is compact it follows that there exists
$y\in \bigcap\limits_{{W\ni x;}\atop W \text{open}}F_W .$ Since $$\bigcap\limits_{{W\ni x;}\atop
W\text{open}}F_W\subseteq K\backslash (\bigcup\limits_{j=1}^kU_j)$$ it follows that
$\sigma(y)\neq x$. Choose an open set $W_x\ni x$ such that $\sigma(y)\notin
\overline{W_x}$. Then $y\notin F_{W_x}$, which is an absurd. This proves the claim.\fimaf

Let $V_0\ni x$ be an open subset according to the claim and define
$$V=V_0\bigcap(\bigcap\limits_{j=1}^k\sigma(U_j)).$$ Let $(y_i)_i$
an net such that $y_i\rightarrow x$. We may suppose that $(y_i)_i\subseteq V$, and so
$\sigma^{-1}(y_i) = \{y_{1,i},... ,y_{k,i}\}$ where $y_{j,i}\in U_j$. How
$\sigma_{|_{U_j}}$ is a homeomorphism we have that
$y_{j,i}\stackrel{i\rightarrow\infty}{\longrightarrow}x_j$ for each $j$, and so
$$L(f)(y_i)=\sum\limits_{z\in U \atop
\sigma(z)=y_i}f(z)=\sum\limits_{j=1}^kf(y_{j,i})\stackrel{i\rightarrow
\infty}{\longrightarrow} \sum\limits_{j=1}^kf(x_j)=\sum\limits_{y\in U\atop
\sigma(y)=x}f(y)=L(f)(x).$$This shows that $L(f)$ is continuous on the points of
$\sigma(K)$, and the lemma is poved. \fim

Now we are in the situation where $C_c(U)$ is an idempotent self-adjoint ideal of
$C_0(U)$, which is an ideal of $C(X)$, and by the previous lemma, $L:C_c(U)\rightarrow
C(X)$ is a function. Moreover, composing $\al$ with the *-isomorphism
$C^b(U)\ni g\mapsto (L_g,R_g)\in M(C_0(U))$ we obtain the partial endomorphism $\widetilde{\al}:C(X)\rightarrow
M(C_0(U))$. It is easy to verify that $(C(X),\widetilde{\al},L)$ is a $C^*$-dynamical
system.

Since $\til{\al}$ is essentialy given by $\al$ we will use the notation
$(C(X),\al,L)$ to us refer to the $C^*$-dynamical system $(C(X),\til{\al},L)$. Moreover,
since $g\til{\al}(f)=g\al(f)$ for each $g\in C_c(U)$ and $f\in C(X)$, no more references
will be made to $\widetilde{\al}$. So we have the Toeplitz algebra $\tp{C(X)}$ and the
crossed product by a partial endomorphism $\prcr{C(X)}$. From now on we will denote
$\tp{C(X)}$ by $\tp{X}$ and $\prcr{C(X)}$ by $\prcr{X}$.

\subsection{Basic results}

Here we will prove some basic results about the crossed product by a partial
endomorphism $\prcr{X}$.

\begin{lema}\label{ignorm} Given $f\in C_c(U)$, we have that: \nl a)$\m{f}=0$
if and only if $f=0$.\nl b)if $\sigma_{|\supp(f)}$ is a homeomorphism then $\|
f\|_\infty=\|\m{f}\|$.
\end{lema}

\demo a) Given $f\in C_c(U)$ and $x\in U$ such that $f(x)\neq 0$ then
$$L(f^*f)(\sigma(x))=\sum\limits_{\sigma(y)=\sigma(x)}f^*(y)f(y)=\sum\limits_{y\neq
x\atop \sigma(y)=\sigma(x)}|f(y)|^2 +|f(x)|^2 > 0.$$ This shows that $L$ is faithful, and
so $\m{f}=0$ if and only if $f=0$. \nl b) Since $\| \m{f}\|^2=\| L(f^*f)\|_\infty$ it
suffices to show that $\| L(f^*f)\|_\infty =\| f\|_\infty^2$. For this note that
$$L(f^*f)(x)=\left\{\begin{array}{ll} |f(\sigma^{-1}(x))|^2 & \text{ if }x\in
\sigma(\supp(f))\\ 0 & \text{ otherwise}
\end{array}\right..$$ Then $\| L(f^*f)\|_\infty\leq \| f\|_\infty^2$. On the other hand,
choose $x\in U$ such that $|f(x)|=\| f\|_\infty$, and note that
$L(f^*f)(\sigma(x))=(f^*f)(x)$, which means that
$\| L(f^*f)\|_\infty\geq \| f\|_\infty^2$. \fim

Consider the *-homomorphism $\varphi:C(X)\rightarrow L(M)$ given by the left product of $A$ by $M$.
Note that $f\in \text{ker}(\varphi)$ if and only if $fm=0$ for each $m\in M$, which
occurs if and only if $\m{fg}=f\m{g}=0$ for each $g\in C_c(U)$. By a) of the previuos
lemma $\m{fg}=0$ if and only if $fg=0$. Therefore $f\in \text{ker}(\varphi)$ if and only
if $fg=0$ for every $g\in C_c(U)$ and so $fg=0$ for all $g\in C_0(U)$. So, given $g\in
C_0(U)$ it follows that $fg=0$ for every $f\in\text{ker}(\varphi)$ and so $f\in
\text{ker}(\varphi)^\bot$. This means that $C_0(U)\subseteq
\text{ker}(\varphi)^\bot$.

\begin{lema}\label{l1}
a) If $f,g\in C_c(U)$ and $\sigma_{|_{\supp(f)\cup \supp(g)}}$ is a homeomorphism then
$(fg^*,\m{f}\m{g}^*)$ is a redundancy of $\tp{X}$ and $fg^*=\m{f}\m{g}^*$ in
$\prcr{X}$.\nl b) $C_0(U)\subseteq \varphi^{-1}(K(M)).$ \nl c) $C_0(U)\subseteq
I_0\,\,(=\varphi^{-1}(K(M))\cap \text{ker}(\varphi)^\bot)$ \nl d) $C_0(U)\subseteq K_1$.
\end{lema}

\demo a) Let $f,g\in C_c(U)$ such that $\sigma_{|_{\supp(f)\cup\supp(g)}}$ is a
homeomorphism and $h\in C_c(U)$. Notice that $\m{f}\m{g}^*\m{h}=\me{f\al(L(g^*h))}$. Since $\sigma_{|_{\supp(f)\cup
\supp(g)}}$ is a homeomorphism, for
each element $x\in \text{supp}(f)$ we have that $f(x)\sum\limits_{y\in
U\atop\sigma(y)=\sigma(x)}(g^*h)(y)=f(x)g(x)^*h(x)$. Therefore for these $x$,
$$f\al(L(g^*h))(x)=f(x)L(g^*h)(\sigma(x))=f(x)\sum\limits_{y\in
U\atop\sigma(y)=\sigma(x)}(g^*h)(y)=f(x)g^*(x)h(x)=(fg^*h)(x).$$ If $x\notin
\text{supp}(f)$ then $(f\al(L(g^*h)))(x)=0=(fg^*h)(x)$. Therefore $f\al(L(g^*h))=fg^*h$.
Then $\m{f}\m{g}^*\m{h}=\me{f\al(L(g^*h))}=\m{fg^*h}=fg^*\m{h}$ for every $h\in
C_c(U)$, from where $\m{f}\m{g}^*m=fg^*m$ for all $m\in M$. It follows that
$(fg^*,\m{f}\m{g}^*)$ is a redundancy. Since $fg^*\in C_0(U)\subseteq
\text{ker}(\varphi)^\bot$ we have that $fg^*=\m{f}\m{g}^*$ in $\prcr{X}$.\nl b) It is
enough to show that $C_c(U)\subseteq K(M)$. Let $f\in C_c(U)$, choose a cover $V_1,\cdots
,V_n$ of $\supp(f)$ such that $\sigma_{|_{V_i}}$ is a homeomorphism. Let $\xi_i''$ be a
partition of unity relative to this cover. Define $\xi_i=f\sqrt{\xi_i''}$ and
$\xi_i'=\sqrt{\xi_i''}$. Then $f=\sum_{i=1}^n\xi_i\xi_i'^*$. By a),
$(\xi_i\xi_i'^*,\m{\xi_i}\m{\xi_i'}^*)$ is a redundancy from where $(f,k)$ is a
redundancy where $k=\sum\limits_{i=1}^n\m{\xi_i}\m{\xi_i'}^*\in \widehat{K_1}$. In this
way $fm=km$ for all $m\in M$ and so $\varphi(f)(m)=fm=km=S(k)(m)$ for every $m\in M$,
where $S$ is the *-isomorphism of \ref{pe1}. It follows that
$\varphi(f)=S(k)$ and so $f\in \varphi^{-1}(K(M))$. Therefore
$C_c(U)\subseteq \varphi^{-1}(K(M))$.\nl c) Follows by b) and by the fact that
$C_0(U)\subseteq \text{ker}(\varphi)^\bot$.
 \nl d) Given $f\in C_c(U)$, by b) it
follows that $(f,k)$ is a redundancy for some $k\in \widehat{K_1}$. Since $f\in
C_0(U)\subseteq I_0$ it follows that $f=q(k)\in K_1$. So $C_c(U)\subseteq
K_1$ from where $C_0(U)\subseteq K_1$. \fim

The following lemma will be used several times in this work.

\begin{lema}\label{multired}
If $(k_0,k_1,\cdots ,k_n)\in C(X)\times K_1\times \cdots \times K_n$ such that
$g\sum\limits_{i=0}^nk_i=0$ for each $g\in C_0(U)$ then:\nl a)
${k_0}_{|_{\partial(U)}}=0$, $k_0=f_1+f_2$ where $f_1\in C_0(U)$ and $f_2\in
C_0(X\setminus \ovl{U})$.\nl b) $\sum\limits_{i=0}^nk_i=f_2$.
\end{lema}

\demo Let $\varepsilon > 0$ be fixed. For every $i\geq 1$ choose
$$k_i'=\sum\limits_{j=1}^{N_i}m_{j,1}^i\cdots m_{j,i}^i(l_{j,1}^i)^*\cdots (l _
{ j , i }^i)^ * \in K_i$$ such that $m_{j,k}^i=\m{f_{j,k}^i}$ with $f_{j,k}^i\in C_c(U)$
and
$\| k_i-k_i'\|\leq\frac{\varepsilon}{n}$. Define
$k_\varepsilon=k_1'+\cdots +k_n'$ and
$K_\varepsilon=\bigcup\limits_{i,j,k}\supp(f_{j,k}^i)\subseteq U$ which is
compact. Given $x\in U\setminus K_\varepsilon$ take $f\in C_0(U)$ such that $f(x)=1$,
$0\leq f\leq 1$ and $f_{|_{K_\varepsilon}}=0$. Then $fk_\varepsilon=0$ by the choice
of $f$ and $fk_0=-f\sum\limits_{i=1}^nk_i$ by hypothesis. It follows that
$$\| fk_0\|=\|-f\sum\limits_{i=1}^nk_i + f_xk_\varepsilon\|= \|
f(-\sum\limits_{i=1}^nk_i + k_\varepsilon)\|=\|f\sum\limits_{i=1}^n(k_i'-k_i)\|\leq
\varepsilon$$ from where $|k_0(x)|\leq \varepsilon$. In this way we have showed that
$|k_0(x)|\leq \varepsilon$ for all $x\in U\setminus K_\varepsilon.$ Given $y\in
\partial(U)$, take a net $(x_l)_l\subseteq U$ such that $x_l\rightarrow y$.
Since $y\notin K_\varepsilon$ and $U\setminus K_\varepsilon$ is open we may suppose
$(x_l)_l\subseteq U\setminus K_\varepsilon$ from where $|k_0(x_l)|\leq \varepsilon$ for
each $l$. By continuity of $k_0$, $|k_0(y)|\leq \varepsilon$. This shows (taking
$\varepsilon$ sufficiently small) that ${k_0}_{|_{\partial(U)}}=0$. Defining $f_1=k_01_U$
and $f_2=k_01_{U^c}$, we obtain a). \nl We will show b). For each $\varepsilon
> 0$ choose $g_\varepsilon\in C_0(U)$ such that $0\leq g\leq 1$ and $g_{|_{K_\varepsilon}}=1$.
Define $h_\varepsilon=g_\varepsilon k_0$. So we obtain a set of functions
$(h_\varepsilon)_\varepsilon\subseteq C_0(U)$.

\af{}{$\lim\limits_{\varepsilon\rightarrow 0}h_\varepsilon=f_1$}

For each $\varepsilon$, given $x\in X$,
$$|(h_\varepsilon-f_1)(x)|=|(g_\varepsilon-1_U)(x)k_0(x)|=\left
\{\begin{array}{cc}|g_\varepsilon(x)-1||k_0(x)| & \text{if }x\in U\setminus K_\varepsilon
\\ 0 & x\in K_\varepsilon\cup U^c
\end{array}\right .$$ For $x\in U\setminus K_\varepsilon$ it holds that $|k_0(x)|\leq \varepsilon$
and so for such elements $|g_\varepsilon(x)-1||k_0(x)|\leq 2\varepsilon.$ So
$\|h_\varepsilon -f_1\|\leq 2\varepsilon$. This shows the claim.\fimaf

Notice that $g_\varepsilon k_\varepsilon=k_\varepsilon$ and $h_\varepsilon=g_\varepsilon
k_0=-g_\varepsilon(k_1+\cdots +k_n)$ because $g_\varepsilon\in C_0(U)$. Then 

$$h_\varepsilon+(k_1+\cdots
+k_n)=h_\varepsilon+k_\varepsilon-k_\varepsilon+(k_1+\cdots +k_n)=$$ $$=-g_\varepsilon(k_1+\cdots
+k_n)+k_\varepsilon-k_\varepsilon+(k_1+\cdots +k_n)=-g_\varepsilon(k_1+\cdots +k_n
-k_\varepsilon)-k_\varepsilon+(k_1+\cdots +k_n),$$ and so $$\| h_\varepsilon+(k_1+\cdots
+k_n)\|=\| g_\varepsilon(-(k_1+\cdots +k_n)+k_\varepsilon)+ ((k_1+\cdots
+k_n)-k_\varepsilon)\|\leq $$ $$\leq\| g_\varepsilon(-(k_1+\cdots +k_n)+k_\varepsilon)\|+
\|(k_1+\cdots +k_n)-k_\varepsilon\|\leq 2\varepsilon.$$ This shows that
$\lim\limits_{\varepsilon\rightarrow 0}h_\varepsilon=-(k_1+\cdots +k_n)$. By the claim
$\lim\limits_{\varepsilon\rightarrow 0}h_\varepsilon=f_1$, and so
$f_1=-(k_1+\cdots +k_n)$. Then
$\sum_{i=0}^nk_i=f_1+f_2+k_1+\cdots+k_n=f_2$, proving b). \fim

\begin{corolario}\label{ig.id}
$K_1\cap C(X)=C_0(U)$
\end{corolario}

\demo Let $r\in K_1\cap C(X)$. Then $r=f=k$ where $f\in C(X)$ and $k\in K_1$. Then
$f-k=0$ and so $g(f-k)=0$ for all $g\in C_0(U)$, and so by $\ref{multired}$,
$f=f_1+f_2$ with $f_1\in C_0(U)$, ${f_2}\in C_0(X\setminus \ovl{U})$ and $f-k=f_2$. How
$f-k=0$ it follows that $f_2=0$. Therefore $f=f_1$, which means that $r=f_1\in
C_0(U)$. In this way $K_1\cap C(X)\subseteq C_0(U)$. The other inclusion is the lema
$\ref{l1}$ d). \fim

In the construction of $\prcr{X}$ we have considered the ideal
$I_0=\varphi^{-1}(K(M))\cap \text{ker}(\varphi)^\bot$. The previous corollary
allows us to identify this ideal.

\begin{corolario}\label{identificacao do ideal}
$I_0=C_0(U)$
\end{corolario}

\demo Given $f\in I_0$ then $\varphi(f)=k\in K(M)$. Choose $k'\in
\widehat{K_1}$ such that $S(k')=k$ where $S$ is the *-isomorphism of \ref{pe1}.
Then $fm=\varphi(f)(m)=k(m)=S(k')(m)=k'm$ for all
$m\in M$. Therefore $(f,k')$ is a redundancy. Since $f\in I_0$ it
follows that $f=q(k')\in K_1$ in $\prcr{X}$. By the previous corollary we have
that $f\in C_0(U)$. So $I_0\subseteq C_0(U)$. The reverse inclusion follows by \ref{l1} c).\fim

Recall that $K$ is the fixed point algebra of the gauge action and that
$K=\ovl{\bigcup\limits_{n\in
\N}L_n}$ where $L_n=C(X)+K_1+\cdots +K_n$ for $n\geq 1$ and $L_0=C(X)$.

\begin{proposicao}\label{int ideal com C(X)}
Every ideal of $\prcr{X}$ which has nonzero intersection with $K$ has nonzero
intersection with $C(X)$.
\end{proposicao}

\demo Let $I$ be an ideal of $\prcr{X}$ such that $I\cap K\neq 0$. By [\ref{davidson}:
III.4.1] there exists $n\in \N$ such that $I\cap L_n\neq 0$. Let $n_0=min\{n\in \N:I\cap
L_n\neq 0\}$ and choose $0\neq k\in I\cap L_{n_0}$. Suppose $n_0\neq 0$. Supposing
$m^*kk^*l=0$ for all $m,l\in M$ we have that $m^*k=0$ for all $m\in M$. So $K_1k=0$ and
by the fact that $C_0(U)\subseteq K_1$ it follows that $fk=0$ for all $f\in C_0(U)$. By
\ref{multired}, $k\in C(X)=L_0$, which is a contradiction because we are supposing $n_0\neq 0$.
So there exists $m,l\in M$ such that $m^*kk^*l\neq 0$. Notice that $m^*kk^*l\in I\cap
L_{n_0-1}$ which again is an absurd because $n_0=\text{min}\{n\in \N:I\cap
L_n\neq 0\}$. Therefore $n_0=0$, that is, $k\in L_0=C(X)$. \fim

By this proposition and by $\ref{idintk}$ follows the corollary:

\begin{corolario}\label{int.de.id.gauge.inv.com.C(X)}
If $0\neq I$ is a gauge-invariant ideal of $\prcr{X}$ then $I\cap C(X)\neq 0$.
\end{corolario}

\subsection{The Cuntz-Krieger algebra for infinite matrices}

We show that that the Cuntz-Krieger algebra for infinite
matrices, introduced in \cite{exelmatinf}, is an example of crossed product by
partial endomorphism. We begin by presenting a short summary of the
construction of this algebra.

Ler $G$ be a set and $A=A(i,j)_{i,j\in G}$ a matrix where each $A(i,j)\in \{0,1\}$.
Define the universal $C^*$-algebra $\m{O_A}$ generated by a set of partial isometries
$\{S_x\}_{x\in G}$ with the following relations:

\begin{enumerate}
\item $S_i^*S_i$ and $S_j^* S_j$ commute, \item $S_i^*S_j=0$ for all $i\neq j$, \item
$S_i^*S_iS_j=A(i,j)S_j$, \item $\prod\limits_{x\in X}S_x^*S_x\prod\limits_{y\in
Y}(1-S_y^*S_y)=\sum\limits_{j\in G}A(X,Y,j)S_jS_j^*$, whenever $X,Y$ are finite subsets of $G$ such that
$A(X,Y,j):=\prod\limits_{x\in X}A(x,j)(1-\prod\limits_{y\in Y}A(y,j))\neq 0$ only for finitely many $j\in G$.
\end{enumerate}

The Cuntz-Krieger algebra for infinite matrices was defined in \cite{exelmatinf} as
the sub-algebra $O_A$ of $\m{O_A}$ generated by the $S_x$.

Let $\F$ be the free group generated by $G$ and let $\{0,1\}^\F$ be the topological space
(with the product topology), which can also be seen as the set of the subsets of $\F$.
In $\{0,1\}^\F$ consider the set $\Omega_e=\{\xi\subseteq \F;e\in\xi\}$, which is
compact. For each $t\in \F$ define $\Delta_t'=\{\xi\in \Omega_e;t\in \xi\}$, which is an
clopen subset. Denoting by $1_t$ the characteristic function of $\Delta_t'$
consider the set $R_A\subseteq C(\Omega_e)$ formed by the following functions:

\begin{enumerate}
\item $1_x1_y$ for all $x\neq y$, $x,y\in G$, \item $1_{x^{-1}}1_y-A(x,y)1_y$ for all
$x,y\in G$, \item $1_{ts}1_t-1_{ts}$ for $t,s\in \F$ such that $|ts|=|t|+|s|$, (where $|s|$
is the number of generators of the reduced form of $s$), \item $\prod\limits_{x\in
X}1_{x^{-1}}\prod\limits_{y\in Y}(1-1_{y^{-1}})-\sum\limits_{j\in G}A(X,Y,j)1_j$ where
$X,Y$ are finite subsets of $G$ such that $A(X,Y,j)\neq 0$ only for finitely many $j\in G$.
\end{enumerate}

In $\Omega_e$ consider the closed set $\m{\Omega_A}=\{\xi\in \Omega_e;f(t^{-1}\xi)=0\,\,
\forall \,\,t\in \xi, f\in R_A \}.$ In [\ref{exelmatinf}:7.3] it was showed that
$\m{\Omega_A}$ is the closure in $\Omega_A^\mathcal{T}$ of the set of the elements which
have an infinite stem (see [\ref{exelmatinf}:5.5]), where
$$\Omega_A^\mathcal{T}=\left\{\begin{array}{ll}
  \xi\in\Omega_e: & e\in \xi,\xi \text{ is convex} \\
   & \text{if }t\in\xi \text{ there is at most one }x\in G
   \text{ such that }tx\in\xi \\
   & \text{if } t\in\xi, y\in G \text{ and } ty\in \xi \text{ then }
   tx^{-1}\in\xi\Leftrightarrow A(x,y)=1
\end{array}\right\}$$

The homeomorphisms $h_t:\Delta_{t^{-1}}'\rightarrow \Delta_t'$ given by $h_t(\xi)=t\xi$
induces a partial action \lb$(\{D_t\}_{t\in \F}, \{\theta_t\})$  (see \cite{exelprodcruzpartial} and
\cite{mclanahan}) of $\F$ in
$C(\m{\Omega_A})$ where $D_t=C(\Delta_t)$, $\Delta_t=\Delta_t'\cap\m{\Omega_A}$ and
$\theta_t:D_{t^{-1}}\rightarrow D_t$ is given by $\theta(f)=f\circ h_{t^{-1}}$
and so we may consider the partial crossed product $C(\m{\Omega_A})\rtimes_{\theta}\F$ (see \cite{exelprodcruzpartial} and \cite{mclanahan}).

It was showed in [\ref{exelmatinf}:7.10] that there exists a *-isomorphism
$\Phi:\m{O_A}\rightarrow C(\m{\Omega_A})\rtimes_{\theta}\F$ such that
$\Phi(S_x)=1_x\delta_x$.

Based on these informations we will show that $\m{O_A}$ is an example of
crossed product by a partial endomorphism.

Let $U\subseteq \m{\Omega_A}$, $U=\bigcup\limits_{x\in G}\Delta_x$. By the fact that each $\Delta_x$
is open it follows that $U$ is open. Moreover, $U$ is dense in $\m{\Omega_A}$, because
$U$ contains all the elements of $\Omega_A^\mathcal{T}$ which have an infinite stem, and
these elements form a dense set in $\m{\Omega_A}$. Since each $\xi\in U$ contains a
unique $x\in G$, we may define the continuous function
$\sigma:U\rightarrow \m{\Omega_A}$ given by $\sigma(\xi)=x^{-1}\xi$ where $x$ is the unique
element of $G$ which lies in $\xi$. This function is a local homeomorphism (in fact,
$\sigma_{|_{\Delta_x}}:\Delta_x\rightarrow \Delta_{x^{-1}}$ is a homeomorphism). Defining
$\al:C(\m{\Omega_A})\rightarrow C^b(U)$ by $\al(f)=f\circ\sigma$
and
$L:C_c(U)\rightarrow C(\m{\Omega_A})$ by $L(f)(\xi)=\sum\limits_{\eta\in
U\atop \sigma(\eta)=\xi}f(\eta)$
we have that $(C(\m{\Omega_A}),\al,L)$ is a
$C^*$-dynamical system, and so we obtain the algebra $\prcr{\m{\Omega_A}}$ (see section
\ref{dynamicalsystem}).

The next step is to show that the algebras $\prcr{\m{\Omega_A}}$ and $\m{O_A}$ are
isomorphic.

\begin{lema}\label{p1}
a) $L(1_x)=1_{x^{-1}}$ for each $x\in G$.\nl b)$f1_x\al L(1_x g)=1_xfg $ for each $x\in G$
and $f,g\in C(\m{\Omega}_A)$.
\end{lema}

\demo Both a) and b) follow by  direct calculation. To prove the first part notice that \lb$\sigma^{-1}(\xi)=\{x\xi:x^{-1}\in \xi\}$.
\fim

\begin{proposicao}\label{p3}
There exists an unitary *-homomorphism $\psi:\m{O_A}\rightarrow \prcr{\m{\Omega_A}}$ such that
$\psi(S_x)=\m{1_x}$.
\end{proposicao}

\demo We will show that $\psi$ preserves the relations 1-4 which defines
$\m{O_A}$. The first relation follows by the fact that
$\psi(S_x)^*\psi(S_x)=\m{1_x}^*\m{1_x}\in C(\m{\Omega}_A)$. To verify the second
relation note that $1_x1_y=0$ for $x,y\in G$ and $x\neq y$, from where
$\psi(S_x)*\psi(S_y)=\m{1_x}^*\m{1_y}=L(1_x1_y)=0$. The third relation follows
by $\ref{p1}$ a) and by the fact that $1_{x^{-1}}1_y=A(x,y)1_y$ in $\m{\Omega_A}$.
In fact,
$$\psi(S_x)^*\psi(S_x)\psi(S_y)=\m{1_x}^*\m{1_x}\m{1_y}=L(1_x)\m{1_y}=1_{x^{-1}}
\m{1_y}=\m{1_{x^{-1}}1_y}=A(x,y)\m{1_y}=A(x,y)\psi(S_y).$$  Let us verify
the fourth relation. By $\ref{l1}$ a) $1_x=\m{1_x}\m{1_x}^*$ in
$\prcr{\m{\Omega_A}}$. Therefore, also
$\sum\limits_{i=1}^n1_{x_i}=\sum\limits_{i=1}^n\m{1_{x_i}}\m{1_{x_i}}^*$ in
$\prcr{\m{\Omega_A}}$. Let $X,Y\subseteq G$ finite such that $A(X,Y,x_i)\neq
0$ only for $i=1,\cdots ,n$. Then $\prod\limits_{x\in
X}1_x^{-1}\prod\limits_{y\in Y}(1-1_y^{-1})=\sum\limits_{i=1}^n1_{x_i}$ in
$\m{\Omega}_A$ and so
$$\prod\limits_{x\in X}\psi(S_x)^*\psi(S_x)\prod\limits_{y\in
Y}(1-\psi(S_y)^*\psi(S_y))=\prod\limits_{x\in X}1_x^{-1}\prod\limits_{y\in
Y}(1-1_y^{-1})=\sum\limits_{i=1}^n1_{j_i}=$$
$$=\sum\limits_{i=1}^n\m{1_{x_i}}\m{1_{x_i}}^*=\sum\limits_{i=1}^n\psi(S_{x_i
})\psi(S_{x_i})^*= \sum\limits_{x\in G}A(X,Y,x)\psi(S_x)\psi(S_x)^*.$$ \fim

We will show that the *-homomorphism defined in this proposition is a
*-isomorphism. The following lemma will be useful to show that this
*-homomorphism is surjective.

\begin{lema}\label{lOA}
The $C^*$-algebra $B$ generated by $\m{1_x}:x\in G$ in $\prcr{\m{\Omega_A}}$ contains all the
elements of $C(\Omega_e)$ of the form $1_r:e\neq r\in \F$ and moreover $B$ coincides with
the $C^*$-algebra generated by $M$.
\end{lema}

\demo By $\ref{p1}$ a), $\m{1_x}^*\m{1_x}=1_{x^{-1}}$. Given $\beta=x_1^{-1}\cdots
x_n^{-1}\in \F $ with $x_i\in G$, by induction
$$\m{1_{x_n}}^*\cdots \m{1_{x_1}}^*\m{1_{x_1}}\cdots \m{1_{x_n}}=
\m{1_{x_n}}^*1_{x_{n-1}^{-1}\cdots x_1^{-1}}\m{1_{x_n}}=L(1_{x_n}1_{x_{n-1}^{-1}\cdots
x_1^{-1}})=1_{x_n^{-1}\cdots x_1^{-1}}.$$ If $b=yr^{-1}$ with $r=x_1\cdots x_n$ and
$x_i,y\in G$ then $$\m{1_y}\m{1_{x_n}}^*\cdots \m{1_{x_1}}^*\m{1_{x_1}}\cdots
\m{1_{x_n}}\m{1_y}^*=\m{1_y}1_{r^{-1}}\m{1_y}^*=\me{1_y\al(1_{r^{-1}})}\m{1_y}^*.$$
By \ref{l1} a) $\me{1_y\al(1_{r^{-1}})}\m{1_y}^*=1_y\al(1_{r^{-1}})$, and by direct calculation 
$1_y\al(1_{r^{-1}})=1_{yr^{-1}}$. Therefore $1_{yr^{-1}}\in B$ for all $y\in G$,
$r=x_1\cdots x_n$ with $x_i\in G$. The general case, $\beta=sr^{-1}$, with $s=x_1\cdots
x_n, r=y_1\cdots y_m$ and $x_i,y_i\in G$ follows by induction. If $t\in \F$ and $t$ is
not of the form $\beta=sr^{-1}$ like above, then $1_t=0$ em $\m{\Omega_A}$ by
[\ref{exelmatinf}:5.8]. Therefore $1_t\in B$ for all $e\neq t\in \F$. We will
show that $B$ is the algebra generated by $M$. For each $x\in G$,
$\text{span}\{1_x\prod\limits_{s}1_s\}$ is dense in $D_x$ and by \ref{ignorm} b), since
$\sigma_{|_{\Delta_x}}$ is a homeomorphism, it follows that
$\text{span}\{\me{1_x\prod\limits_{s}1_s}\}$ is dense in $\m{D_x}$. Since
$\me{1_x\prod\limits_{s}1_s}=1_x\prod\limits_{s}1_s\m{1_x}\in B$ we have that
$\m{D_x}\subseteq B$, because $B$ is closed. So $\m{C_c(U)}\subseteq B$ and since $B$ is
closed it follows that $M\subseteq B$. This shows that $B$ contains the algebra generated
by $M$. On the other hand, since $\m{1_x}\in M$ for each $x\in G$, it is clear that the
algebra generated by $M$ contains $B$, and this concludes the proof.\fim

\begin{proposicao}\label{hom.entre.prod.cruz.e.prodcruz.parc}
There exists a *-homomorphism $\phi:\prcr{\m{\Omega_A}}\rightarrow C(\m{\Omega}_A)\rtimes_\theta \F$ such
that $\phi(f)=f\delta_e$ for all $f\in C(X)$ and $\phi(\m{f_x})=f_x\delta_x$ for all $f\in D_x$ and $x\in G$.
\end{proposicao}

\demo Let us define initially a homomorphism from the Toeplitz algebra
$\tp{\m{\Omega_A}}$ to $C(\m{\Omega}_A)\rtimes \F$. Define
$\phi':C(\m{\Omega_A})\rightarrow C(\m{\Omega_A})\rtimes_\theta \F$ by $\phi'(f)=f\delta_e$ and
$\phi'':\m{C_c(U)}\rightarrow C(\m{\Omega_A})\rtimes_\theta \F$ by $\phi''(\m{f_x})=f_x\delta_x$ for $f_x\in
D_x$.  Clearly $\phi'$ is a *-homomorphism. By
\ref{ignorm} a) $\phi''$ is well defined. Moreover $\phi''$ is linear and given $g=\sum
g_x$ and $f=\sum f_x$ in $C_c(U)$, where $f_x,g_x\in D_x$, we have that
$$\phi''(\m{g})^*\phi''(\m{f})=(\sum g_x\delta_x)^*(\sum f_y\delta_y)=(\sum \theta_{x^{-1}}(g_x^*)\delta_{x^{-1}})(\sum
f_y\delta_y)=$$
$$=\sum\limits_{x,y}\theta_{x^{-1}}(g_x^*)\delta_{x^{-1}}
f_y\delta_y=\sum\limits_{x,y}\theta_{x^{-1}}(g_x^*f_y)\delta_{x^{-1}y}=\sum\limits_x\theta_{x^{-1}}(g_x^*f_x)\delta_e.$$

\af{}{$L(g^*f)=\sum \theta_{x^{-1}}(g_x^*f_x)$}

It is enough to show that $L(g_x^*f_x)=\theta_x^{-1}(g_x^*f_x)$ because $g_x^*f_y=0$ for
$x\neq y$. For this notice that if $x^{-1}\notin \xi$ then
$L(g_x^*f_x)(\xi)=0=\theta_x^{-1}(g_x^*f_x)(\xi)$. Moreover, if $x^{-1}\in \xi$ then we have
$L(g_x^*f_x)(\xi)=(g_x^*f_x)(x\xi)=(g_x^*f_x)(h_x(\xi))=\theta^{-1}(g_x^*f_x)(\xi)$. So
the claim is proved.\fimaf

Then
$\sum\limits_x\theta_{x^{-1}}(g_x^*f_x)\delta_e=L(g^*f)\delta_e=\phi'(\langle
\m{g},\m{f}\rangle ),$ and so, $\phi''(\m{g})^*\phi''(\m{f})=\phi'(\langle
\m{g},\m{f}\rangle ).$ Therefore
$$\|\phi''(\m{f})\|^2=\|\phi''(\m{f})^*\phi''(\m{f})\|=\|\phi'(\langle
\m{f},\m{f}\rangle )\|\leq\| \langle \m{f},\m{f}\rangle \|=\|\m{f}\|_M^2$$ from where we
may extend $\phi''$ to $M$. In this way we obtain a function
$$\phi:C(\m{\Omega_A})\cup M\rightarrow C(\m{\Omega}_A)\rtimes_{\theta} \F$$
defined by $\phi(f)=\phi'(f)$ if $f\in C(\m{\Omega_A})$ and $\phi(m)=\phi''(m)$ for $m\in
M$.

\af{}{$\phi$ satisfies the relations which defines $\tp{\m{\Omega_A}}$.}

By density of $\m{C_c(U)}$ in $M$ it suffices to verify if $\phi$ satisfies the relations
for elements of the form $\m{f}=\sum\m{f_x},\m{g}=\sum\m{g_y}\in \m{C_c(U)}$, where $f_x, g_x\in D_x$,  and $h\in C(\m{\Omega_A})$. We
already know that $\phi$ preserves the relations of $C(\m{\Omega_A})$, of $M$ and that
$\phi(\m{f})^*\phi(\m{g})=\phi(\langle \m{f},\m{g}\rangle )$. Moreover,
$$\phi(h)\phi(\m{f})=h\delta_e\sum f_x\delta_x=\sum
hf_x\delta_x=\phi(\m{hf})=\phi(h\m{f})$$ and
$$\phi(\m{f})\phi(h)=(\sum f_x\delta_x)h\delta_e=\sum
\theta_x(\theta_x^{-1}(f_x)h)\delta_x=\sum f_x\al(h)\delta_x=\phi(\m{f\al(h)})=\phi(\m{f}h).$$ This proves the claim. \fimaf

Se we may
extend $\phi$ to $\tp{\m{\Omega_A}}$. We will show that if $(a,k)$ is a
redundancy then $\phi(a)=\phi(k)$. For each $f_x\in D_x$,
$\phi(\m{f_x}\m{1_x}^*)=f_x\delta_x1_{x^{-1}}\delta_{x^{-1}}=f_x\delta_e=\phi(f_x)$
and so if $f=\sum\limits_x f_x$ with $f_x\in D_x$ then
$\phi(f)=\sum\limits_x\phi(\m{f_x}\m{1_x}^*)$. Given a redundancy $(f,k)$ with $f\in
I_0$, and so $f\in C_0(U)$ by \ref{identificacao do ideal}, choose $(f_n)_n\subseteq
C_c(U)$ such that $f_n\rightarrow f$, and $(k_n)_n\subseteq \widehat{K_1}$ such that
$k_n\rightarrow k$ and $k_n=\sum\limits_{i=1}^{t_n}m_{i,n}r_{i,n}^*$ with
$m_{i,n},r_{i,n}\in M$. Since $f_n\in C_c(U)$ for each $n$, we have that
$f_n=\sum\limits_{i=1}^{l_n} f_{x_{i,n}}$ and so $\phi(f_n)=\sum\limits_{i=1}^{l_n}
\phi(\m{f_{x_{i,n}}}\m{1_{x_{i,n}}}^*)$. Then

$$\phi(f-k)\phi(f-k)^*=\lim\limits_n\phi(f-k)(\phi(f_n)^*-\phi(k_n^*))=$$ $$=\lim\limits_n\phi(f-k)(\phi(\sum\limits
\m{1_{x_{i,n}}}\m{f_{x_{i,n}}}^*)-\phi(\sum\limits_{i=1}^{n_i}r_{i,n}m_{i,n}^*))=$$ $$=\lim\limits_n\phi((f-k)(\sum\limits
\m{1_{x_i}}\m{f_{x_i}}^*-\sum\limits_{i=1}^{n_i}r_{i,n}m_{i,n}^*))=0.$$
The last equality follows by the fact that $(f-k)m=0$ for each $m\in M$, because $(f,k)$ is a redundancy.
This shows that
$\phi(f)=\phi(k)$. \fim

\begin{proposicao}\label{teoisock}
The *-homomorphism $\psi:\m{O_A}\rightarrow \prcr{\m{\Omega_A}}$ defined in $\ref{p3}$
is a \lb*-isomorphism.
\end{proposicao}

\demo To prove that $\psi$ is surjective it is anough to prove that $M\cup
C(\m{\Omega_A})\subseteq \text{Im}(\psi)$. By the lemma $\ref{lOA}$, $M\subseteq
\text{Im}(\psi)$. By the same lemma, the elements of the form $1_r:e\neq r\in \F$ are in
the range of $\psi$ and moreover, $\psi(1)=1=1_e$. The algebra generated by the elements
$\{1_r:r\in F\}$ is self-adjoint, contains the constant functions and separate points,
and so is dense in $C(\m{\Omega_A})$. It follows that $C(\m{\Omega_A})\subseteq
\text{Im}(\psi)$. In order to see that $\psi$ is injective, note that $\Phi^{-1}\phi\psi=Id_{\m{O}_A}$
where $\phi$ is the *-homomorphism of $\ref{hom.entre.prod.cruz.e.prodcruz.parc}$ and
$\Phi$ is the *-isomorphism between $\m{O_A}$ and $C(\m{\Omega_A})\rtimes_\theta \F$ such
that $\Phi(S_x)=1_x\delta_x$. \fim

By this proposition and by $\ref{lOA}$ it follows that the Cuntz-Krieger algebra for
infinite matrices $O_A$ is isomorphic to the algebra $B$, generated by $M$. Note that the
algebra generated by $M$ coincides with the ideal $\langle M\rangle $ of
$\prcr{\m{\Omega_A}}$.

\section{Relationship between the gauge-invariant ideals of $\prcr{X}$ and open sets of
X}

We show in this section a bijection between the gauge-invariant ideals
of $\prcr{X}$ and the $\sigma,\sigma^{-1}$-invariant subsets of $X$. In
particular, we prove that every gauge-invariant ideal of $\prcr{X}$ is
generated by the set $C_0(V)$ for some $V\subseteq X$ whith is
$\sigma,\sigma^{-1}$-invariant.

\begin{definicao}\label{sigma.sigma-1.inv}
a) A set $V\subseteq X$ is $\sigma$-invariant if $\sigma(V\cap U)\subseteq V$.\nl b) A
set $V\subseteq X$ is $\sigma^{-1}$-invariant if $\sigma^{-1}(V) \subseteq V$.\nl c) A
set $V\subseteq X$ is $\sigma,\sigma^{-1}$-invariant if it is $\sigma$-invariant and
$\sigma^{-1}$-invariant.
\end{definicao}

Let $V\subseteq X$ be an open set. We say that $C_0(V)$ is
$L$-invariant if $L(C_0(V)\cap C_c(U))\subseteq C_0(V)$.

\begin{proposicao}\label{sigma.inv.ss.L.inv.}
a)A open set $V\subseteq X$ is $\sigma$-invariant if and oly if $C_0(V)$ is
$L$-invariant.\nl b)A open set $V\subseteq X$ is $\sigma^{-1}$-invariant if and only if
$f\al(g)\in C_0(V)$ for all $f\in C_c(U)$ and $g\in C_0(V)$.
\end{proposicao}

\demo a) Suppose $V$ $\sigma$-invariant. Given $f\in
C_0(V)\cap C_c(U)$, choose $x\notin V$. Supposing $y\in \sigma^{-1}(x)\cap V$, we
have $x=\sigma(y)\in V$ because $V$ is $\sigma$-invariant. So there does not exists a such $y$, and therefore 
$L(f)(x)=0$. This shows that $L(f)\in
C_0(V)$. On the other hand, suppose $C_0(V)$ $L$-invariant. Suppose $x\in U\cap V$
and choose $f_x\in C_c(U)\cap C_0(V)$ such that $f_x(x)\neq 0$. Then $L(f_x^*f_x)\in C_0(V)$ and
$L(f_x^*f_x)(\sigma(x))\neq 0$, which shows that
$\sigma(x)\in V$.\nl b) Suppose $V$
$\sigma^{-1}$-invariant. Let $f\in C_c(U)$, $g\in
C_0(V)$ and $x\notin V$. If $x\notin U$, then $f(x)=0$ and so $(f\al(g))(x)=0$. If $x\in U$, since $V$ is
$\sigma^{-1}$-invariant then
$\sigma(x)\notin V$ and therefore $f\al(g)(x)=f(x)g(\sigma(x))=0$. So $f\al(g)\in C_0(V)$. On the other hand, let $x\in
\sigma^{-1}(y)$, $y\in V$. Choose $g\in C_0(V)$ such that $g(y)\neq 0$ and
$f\in C_c(U)$ such that $f(x)\neq 0$. Then, since $f\al(g)\in C_0(V)$ and
$(f\al(g))(x)=f(x)g(y)\neq 0$ it follows that $x\in V$. So $V$ is $\sigma^{-1}$-invariant.
\fim

If $V\subseteq X$ is an open $\sigma, \sigma^{-1}$-invariant set 
then $X'=X\setminus V$ is a compact $\sigma, \sigma^{-1}$-invariant set. Define $U'=U\cap
X'\,\,(=U\setminus V)$ and consider $\sigma':=\sigma_{|_U'}:U'\rightarrow X'$ which is a local homeomorphism. 
Consider the $C^*$-dynamical system $(X',\al',L')$ where $\al'$ and $L'$ are defined as $\al$ and $L$ in the section 
\ref{dynamicalsystem}. Denote by $M'$
the Hilbert module generated by $C_c(U')$, by $\langle C_0(V)\rangle $ the ideal generated by
$C_0(V)$ in $\prcr{X}$ and by $\ovl{\ovl{b}}$ the image of the elements $b\in
\prcr{X}$ by the quotient map of $\prcr{X}$ on $\prcr{X} /\langle C_0(V)\rangle
$.

\begin{teorema}\label{teoiso}
There exists a *-isomorphism $\Psi:\prcr{X}/\langle C_0(V)\rangle \rightarrow
\mathcal{O}(X',\al',L')$ such that $\Psi(\ovl{\ovl{f}})=f_{|_{X'}}$ for each $f\in
C(X)$.
\end{teorema}

\demo Define
$\Psi_1:C(X)\rightarrow C(X')$ by $\Psi_1(f)=f_{|_{X'}}$
which is a *-homomorphism and is
surjective, by Tietze's theorem. Moreover, for every $\m{f}\in \m{C_c(U)}\subseteq M$
define
$\Psi_2(\m{f})=\m{f_{|_{X'}}},$ which is a linear and contractive of $\m{C_c(U)}\subseteq M$
to $M'$ and so we may extend it to $M$. So we may define in an obvious manner
$\Psi_3:C(X)\cup M\rightarrow \mathcal{T}(X',\al',L').$ It is easy to verify that $\Psi_3$ satisfies
the relations that defines $\tp{X}$ and so $\Psi_3$ has an extension to $\tp{X}$,
which will be denoted by $\Psi_3$. We will show that $\Psi_3$ is surjective.
Given $h\in C_c(U')$, choose $g\in C_c(U)$ such that $g_{|_{\supp(h)}}=1$ and
$f\in C(X)$ such that $\Psi_3(f)=h$. Then $fg\in C_c(U)$ and
$\Psi_3(f)\Psi_3(\m{g})=h\m{g_{|_{X'}}}=\m{hg_{|_{X'}}}=\m{h}$. This shows that $\Psi_3(M)$ is
dense in $M'$, and with the fact that $C(X')\subseteq \text{Im}(\Psi_3)$, it
follows that $\Psi_3$ is surjective.

\af{}{If $(f,k)$ is a redundancy of $\tp{X}$ and $f\in I_0$ then
$(\Psi_3(f),\Psi_3(k))$ is a redundancy of $\mathcal{T}(X',\al',L')$ and $\Psi_3(f)\in
I_0'$}.

Let $(f,k)$ be a redundancy of $\tp{X}$ and $f\in I_0$. Then $fm=km$, from where
$\Psi_3(f)\Psi_3(m)=\Psi_3(k)\Psi_3(m)$. Since $\Psi_3(f)\in C(X')$ and $\Psi_3(k)\in
\widehat{K_1'}$ and moreover $\Psi_3(M)$ is dense in $M'$ it follows that
$(\Psi_3(f),\Psi_3(k))$ is a redundancy. Since $f\in I_0$, and $I_0=C_0(U)$ by
\ref{identificacao do ideal}, it follows that $f\in C_0(U)$ and therefore
$\Psi_3(f)=f_{|_{X'}}\in C_0(U')=I_0'$.\fimaf

If $q$ is the quotient map of $\mathcal{T}(X',\al',L')$ on $\mathcal{O}(X',\al',L')$
then the composition $q\circ \Psi_3$ is a \lb*-homomorphism of $\tp{X}$ on
$\mathcal{O}(X',\al',L')$ which by the claim above vanishes on the elements $(a-k)$ for all redundancies $(a,k)$ such
that $a\in I_0$. By passage to the quotient we obtain a *-homomorphism of $\prcr{X}$
to $\mathcal{O}(X',\al',L')$ which will be denoted by $\Psi_0$. Moreover, given $f\in
C_0(V)$ note that $\Psi_0(f)=f_{|_{X'}}=0$, and again passing to the quotient we obtain an
other *-homomorphism of $\prcr{X}/\langle C_0(V)\rangle $ to
$\mathcal{O}(X',\al',L')$, which will be called $\Psi$. It remains to show that
$\Psi$ is injective. Note that $\langle C_0(V)\rangle$ is gauge-invariant. Consider the
gauge action on $\prcr{X}/\langle C_0(V)\rangle$ whose fixed point algebra is
$\ovl{\ovl{K}}=\ovl{\bigcup\limits_{n\in \N}\ovl{\ovl{L_n}}}$ (see paragraph following
\ref{defgaugeinv}) and the gauge action on $\mathcal{O}(X',\al',L')$. Since $\Psi$ is
covariant by these actions, by [\ref{exelprodcruzpartial}: 2.9] it is enough to show that $\Psi$ restricted to
$\ovl{\ovl{K}}$ is injective. For this we will show that $\Psi$ restricted to
$\ovl{\ovl{L_n}}$ is injective for all $n\in \N$.

\af{ 1}{let $\ovl{\ovl{k_0}}+\ovl{\ovl{k_1}}+\cdots + \ovl{\ovl{k_n}}\in
\ovl{\ovl{L_n}}$. If $\phi(\ovl{\ovl{k_0}}+\ovl{\ovl{k_1}}+\cdots + \ovl{\ovl{k_n}})=0$
then $\ovl{\ovl{k_0}}\in \ovl{\ovl{K_1}}$.}

Let $k_i'=\Psi(\ovl{\ovl{k_i}})$ and notice that $k_0'\in C(X')$ and $k_i'\in K_i'$ for $i\geq 1$. Then $k_0'+k_1'+\cdots +k_n'=0$ and
so $g(k_0'+k_1'\cdots +k_n')=0$ for all $g\in C_0(U')$. By $\ref{multired}$ it follows
that $k_0'=f_1+f_2$ where $f_1\in C_0(U')$ and $k_0'+k_1'\cdots +k_n'=f_2$ from where
$f_2=0$. Then $k_0'\in C_0(U')$ and so $k_0\in C_0(U\cup V)$ from where
$\ovl{\ovl{k_0}}\in \ovl{\ovl{C_0(U\cup V)}}=\ovl{\ovl{C_0(U)}}+\ovl{\ovl{C_0(V)}}\subseteq \ovl{\ovl{K_1}}.$ 

\af{ 2}{$\Psi$ restricted to $\ovl{\ovl{C(X)}}$ is faithful, and also $\Psi$ restricted
to $\ovl{\ovl{K_n}}$ is faithful.}

If $f\in C(X)$ and $\Psi(\ovl{\ovl{f}})=0$ then $f\in C_0(V)$ and so $\ovl{\ovl{f}}=0$.
This shows the first part. To prove the second assertion let $\ovl{\ovl{k_n}}\in
\ovl{\ovl{K_n}}$ and suppose that $\Psi(\ovl{\ovl{k_n}})=0$. Then
$\Psi({\ovl{\ovl{M}}^*}^n\ovl{\ovl{k_n}}\,{\ovl{\ovl{M}}\,}^n)=0$ and how
${\ovl{\ovl{M}}^*}^n\ovl{\ovl{k_n}}\,{\ovl{\ovl{M}}\,}^n\subseteq
\ovl{\ovl{C(X)}}$ and $\Psi$ restricted to $\ovl{\ovl{C(X)}}$ is faithful it follows
that
${\ovl{\ovl{M}}^*}^n\ovl{\ovl{k_n}}\,{\ovl{\ovl{M}}\,}^n=0$ from where
$\ovl{\ovl{K_n}}\,\ovl{\ovl{k_n}}\,\ovl{\ovl{K_n}}=0$ and so
$\ovl{\ovl{k_n}}=0$.\fimaf

We will prove now the following claim which will conclude the proof of the theorem.

\af{ 3}{For all $n\in \N$, $\Psi$ restricted to $\ovl{\ovl{L_n}}$ is faithful}

By claim 2 $\Psi$ restricted to $\ovl{\ovl{L_0}}$ is faithful. By induction, suppose that
$\Psi$ restricted to $\ovl{\ovl{L_n}}$ is faithful, take
$\ovl{\ovl{k_0}}+\ovl{\ovl{k_1}}+\cdots +\ovl{\ovl{k_{n+1}}}\in
\ovl{\ovl{L_{n+1}}}$ and suppose that
$\Psi(\ovl{\ovl{k_0}}+\ovl{\ovl{k_1}}+\cdots +\ovl{\ovl{k_{n+1}}})=0.$ Then
$$\Psi(\ovl{\ovl{M}}^*(\ovl{\ovl{k_0}}+\ovl{\ovl{k_1}}+\cdots
+\ovl{\ovl{k_{n+1}}})^*(\ovl{\ovl{k_0}}+\ovl{\ovl{k_1}}+\cdots+
\ovl{\ovl{k_{n+1}}})\ovl{\ovl{M}})=0$$ and by the induction hypothesis,
$$\ovl{\ovl{M}}^*(\ovl{\ovl{k_0}}+\ovl{\ovl{k_1}}+\cdots
+\ovl{\ovl{k_{n+1}}})^*(\ovl{\ovl{k_0}}+\ovl{\ovl{k_1}}+.. .+
\ovl{\ovl{k_{n+1}}})\ovl{\ovl{M}}=0,$$ from where
$(\ovl{\ovl{k_0}}+\ovl{\ovl{k_1}}+\cdots + \ovl{\ovl{k_{n+1}}})\ovl{\ovl{M}}=0$ and
so
$(\ovl{\ovl{k_0}}+\ovl{\ovl{k_1}}+\cdots + \ovl{\ovl{k_{n+1}}})(\ovl{\ovl{K_1}}+\cdots
+\ovl{\ovl{K_{n+1}}})=0.$ By claim 1, $\ovl{\ovl{k_0}}\in \ovl{\ovl{K_1}}$, from where
$\ovl{\ovl{k_0}}+\ovl{\ovl{k_1}}+\cdots +\ovl{\ovl{k_{n+1}}}\in(\ovl{\ovl{K_1}}+\cdots
+\ovl{\ovl{K_{n+1}}})$ and therefore $\ovl{\ovl{k_0}}+\ovl{\ovl{k_1}}+\cdots
+\ovl{\ovl{k_{n+1}}}=0.$ \fim

Given an ideal $I$ in $\prcr{X}$, the set $I\cap C(X)$ is an ideal of $C(X)$ and so it is
of the form $C_0(V)$ for some open set $V\subseteq X$. The following proposition shows a
feature of these open sets.

\begin{proposicao}\label{sigma invariante}
Let $I\unlhd\prcr{X}$ and $V\subseteq X$ the open set such that $I\cap C(X)=C_0(V)$. Then $V$ is a $\sigma,
\sigma^{-1}$-invariant set.
\end{proposicao}

\demo Given $f\in C_c(U)\cap C_0(V)$, take $g\in C_c(U)$ such that $g_{|_{\text{supp}(f)}}=1$. Then
$f\m{g}=\m{f}\in I$ and so $L(f)=\m{g}^*\m{f}\in I\cap C(X)=C_0(V)$. By \ref{sigma.sigma-1.inv} a) it follows
that $V$ is $\sigma$-invariant.
We will show that $V$ is a $\sigma^{-1}$-invariant set. Let
$x$ be an element of $V$ and $y\in \sigma^{-1}(x)$. Choose $f_x\in C_0(V)$ such that $f_x(x)=1$
and $f_y\in C_c(U)$ such that $f_y(y)=1$ and $\sigma_{|_{\supp(f_y)}}$ is a
homeomorphism. Then $\me{f_y\al(f_x)}=\m{f_y}f_x\in I\cap M$ and therefore
$\me{f_y\al(f)}\m{f_y}^*\in I$. By \ref{l1} a),
$f_y\al(f_x)f_y^*=\me{f_y\al(f_x)}\m{f_y}^*$ and so $f_y\al(f_x)f_y^*\in I\cap
C(X)=C_0(V)$. Note that
$$(f_y\al(f_x)f_y^*)(y)=|f_y|^2(y)f_x(\sigma(y))=|f_y(y)|^2f_x(x)=1,$$ which
shows that $y\in V$.\fim

This proposition shows that there exists a map
$$\Phi:\{\text{ideals of } \prcr{X}\}\rightarrow \{\text{open
}\sigma,\sigma^{-1}-\text{invariant sets of } X\}$$ given by $\Phi(I)=V$ where $V$ is
the open set of $X$ such that $I\cap C(X)=C_0(V)$. The following proposition shows that
$\Phi$ is surjective. To prove this proposition we need some lemmas.

\begin{lema}\label{pteo1}
Let $V$ a $\sigma$-invariant set and $f_1,\cdots ,f_n,g_1,\cdots ,g_n\in C_c(U)$ such that $f_i\in
C_0(V)$ or $g_i\in C_0(V)$ for some $i$. Then $\m{f_n}^*\cdots \m{f_1}^*\m{g_1}\cdots
\m{g_n}\in C_0(V)$.
\end{lema}

\demo Suppose $f_i\in C_0(V)$ and define $h_{j}=\m{f_j}^*\cdots \m{f_1}^*\m{g_1}\cdots
\m{g_j}$ for $j\geq 1$ and $h_0=1$. Since $h_j\in C(X)$ for each $j$ it follows that $f_i^*h_{i-1}g_i\in
C_0(V)$. By \ref{sigma.inv.ss.L.inv.} $C_0(V)$ is $L$-invariant, and so
$h_i=\m{f_i}^*h_{i-1}\m{g_i}=L(f_i^*h_{i-1}g_i)\in C_0(V)$. By induction it may be showed that $h_n\in C_0(V)$. If $g_i\in C_0(V)$ the
proof is analogous. \fim

To show that the map $\Phi$ is surjective we will show that if $V$ is an open
$\sigma,\sigma^{-1}$-invariant set then $\langle C_0(V)\rangle \cap C(X)=C_0(V)$. The
following arguments are a preparation to prove this fact. Given $f\in \langle
C_0(V)\rangle \cap C(X)$ and $\varepsilon > 0$ then there are $a_i,b_i\in \prcr{X}$,
$h_i\in C_0(V)$ such that $$\| f-\sum\limits_{i=1}^Na_ih_ib_i\|\leq\varepsilon$$ where
each $a_i$ is of the form $a_i=m_1\cdots m_{r_i}n_1^*\cdots n_{s_i}^*$ or $a_i\in C(X)$ and each
$b_i$ is of the form $b_i=p_1\cdots p_{t_i}q_1^*\cdots q_{l_i}^*$ or $b_i\in C(X)$. Moreover we
may suppose that $m_j=\m{z_j}$, $n_j=\m{w_j}$, $p_j=\m{u_j}$, $q_j=\m{v_j}$ for each
$m_j, n_j,p_j,$ and $q_j$. Considering the conditional expectation $E$ induced by the
gauge action and that $$\| f-\sum\limits_{i=1}^NE(a_ih_ib_i)\|=\|
E(f-\sum\limits_{i=1}^Na_ih_ib_i)\|\leq \varepsilon,$$ we may suppose that $r_i+t_i=s_i+l_i$,
because
$$E(a_ih_ib_i)=\left\{\begin{array}{cl} a_ih_ib_i & \text{if }r_i+t_i=s_i+l_i\\ 0 & \text{otherwise}
\end{array}\right ..$$

\begin{lema}\label{pteo2}
Let $V$ be an open $\sigma,\sigma^{-1}$-invariant set. Then for each $i$ we have that $a_ih_ib_i\in
C_0(V)$ or $a_ih_ib_i=\m{f_1}\cdots \m{f_n}\m{g_n}^*\cdots \m{g_1}^*$ where $f_j\in
C_0(V)$ for some $j$ or $g_j\in C_0(V)$ for some $j$.
\end{lema}

\demo Recall that $a_i=\m{z_1}\cdots \m{z_{r_i}}\m{w_1}^*\cdots \m{w_{s_i}}^*$ or $a_i\in C(X)$,
$b_i=\m{u_1}\cdots \m{u_{t_i}}\m{v_1}^*\cdots \m{v_{l_i}}^*$ or $b_i\in C(X)$ and $r_i+t_i=s_i+l_i$.

Suppose $s_i\leq t_i$. By \ref{pteo1} $w=\m{w_1}^*\cdots \m{w_{s_i}}^*h_i\m{u_1}\cdots \m{u_{s_i}}\in
C_0(V)$ (if $s_i=0$ then $w=h_i$). If $t_i\neq s_i$ then write
$a_ih_ib_i=\m{z_1}\cdots \m{z_{r_i}}\m{wu_{{s_i}+1}}\cdots \m{u_{t_i}}\m{v_1}^*\cdots\m{v_{l_i}}^*,$
and note that $wu_{{s_i}+1}\in C_0(V)$ and therefore $a_ih_ib_i$ is in the desired form. If
$t_i=s_i$ then $r_i=l_i$. If $r_i=0$ (and so $l_i=0$) then $a_ih_ib_i=w\in C_0(V)$. If $r_i\neq 0$
write
$a_ih_ib_i=\m{z_1}\cdots \m{z_{r_i}\al(w)}\m{u_{{s_i}+1}}\cdots
\m{u_{t_i}}\m{v_1}^*\cdots\m{v_{l_i}}^*,$ and in this case $z_{r_i}\al(w)\in C_0(V)$ by the fact
that $V$ is $\sigma^{-1}$-invariant, and so $a_ih_ib_i$ is in the desired form.

Supposing $s_i>t_i$ consider the element $(a_ih_ib_i)^*$, which is in the desired form of the
lemma by the proof above, and therefore $a_ih_ib_i$ is also in the desired form. \fim

The following lemma is only a summary from \ref{pteo1} to \ref{pteo2}.

\begin{lema}\label{res.pteo1.pteo2}
If $V$ is an open $\sigma,\sigma^{-1}$-invariant set then given $f\in \langle
C_0(V)\rangle \cap C(X)$ and $\varepsilon > 0$, there exists $d_0\in C_0(V)$ and
$d_i=\m{f_1^i}\cdots \m{f_{n_i}^i}\m{g_{n_i}^i}^*\cdots \m{g_1^i}^*$, with $f_j^i\in
C_0(V)$ or $g_j^i\in C_0(V)$ for some $j$, $i=1,\cdots ,N$, such that
$\|f-(d_0+\sum\limits_{i=1}^N d_i)\|\leq \varepsilon.$
\end{lema}

Now we prove the proposition which shows that the map $\Phi$ is surjective.

\begin{proposicao}\label{prop.CoV}
If $V\subseteq X$ is $\sigma,\sigma^{-1}$-invariant then $\langle
C_0(V)\rangle \cap C(X)=C_0(V)$.
\end{proposicao}

\demo It is clear that $C_0(V)\subseteq \langle C_0(V)\rangle \cap
C(X)$. To show that $\langle C_0(V)\rangle \cap C(X)\subseteq
C_0(V)$ we will show that given $f\in \langle C_0(V)\rangle \cap
C(X)$, for every $\varepsilon > 0$ it holds that $|f(x)|\leq \varepsilon$ for
each $x\notin V$.

Given $f\in \langle C_0(V)\rangle \cap C(X)$ and $\varepsilon > 0$, by
\ref{res.pteo1.pteo2} we may consider $\|
f-(d_0+\sum\limits_{i=1}^Nd_i)\|\leq \varepsilon$ with $d_0\in C_0(V)$,
$d_i=\m{f_1^i}\cdots \m{f_{n_i}^i}\m{g_{n_i}^i}^*\cdots \m{g_1^i}^*$ where
$f_j^i\in C_0(V)$ for some $j$ or $g_j^i\in C_0(V)$ for some $j$. Define
$K=\bigcup_{i=1}^N\bigcup_{j=1}^{n_i}(\supp(f_j^i)\cup \supp(g_j^i))$ which
is a compact subset of $U$.

 \af{ 1}{If $x\notin V$ and $x\notin U$ then
$|f(x)|\leq \varepsilon$}

If $x\notin U$, choose $h\in C(X)$, $0\leq h\leq 1$, such that
$h(x)=1$ e $h_{|_K}=0$. Then $hd_i=0$ for $i\geq 1$ and so
$\| h(f-d_0)\| =\| h(f-d_0+\sum\limits_{i=1}^Nd_i)
\|\leq\varepsilon$ from where
$|f(x)-d_0(x)|=|(h(f-d_0))(x)|\leq\varepsilon$. Since $x\notin V$
it follows that $d_0(x)=0$ and therefore $|f(x)|\leq\varepsilon.$\fimaf

Now we study the case $x\notin V$ and $x\in U$. Let
$N_0=max\{n_1,\cdots ,n_N\}$. Supposing $N_0=0$, that is, $d_i=0$ fore each
$i\geq 1$, we have that $|f(x)|=|f(x)-d_0(x)|\leq \varepsilon$. Suppose
therefore that $N_0\geq 1$. Let us analyse the case
$\sigma^{N_0-1}(x)\in U$. Define $x_j=\sigma^j(x)$ for $j\in \{0,\cdots
,N_0\}$. For each $j\in \{0,\cdots,N_0-1\}$ take $h_j\in C_c(U)$ such that
$h_j(x_j)=1$, $0\leq h_j\leq 1$ and $\sigma_{\supp(h_j)}$ is a homeomorphism.

\af{ 2}{For each $i\in\{0,\cdots, N\}$,
$h_i'=\m{h_{N_0-1}}^*\cdots \m{h_0}^*d_i\m{h_0}\cdots \m{h_{N_0-1}}\in C_0(V)$.}

For $i\geq 1$, since $f_j^i\in C_0(V)$ or $g_j^i\in C_0(V)$ for some $j$, by
\ref{pteo1} we have that $$u=\m{h_{n_i-1}}^*\cdots \m{h_0}^*\m{f_1^i}\cdots
\m{f_{n_i}^i}\in C_0(V) \text{ or }
v=\m{g_{n_i}^i}^*\cdots \m{g_1^i}^*\m{h_0}\cdots \m{h_{n_i-1}}\in
C_0(V).$$ Then $uv\in C_0(V)$ and again by \ref{pteo1} it follows that
$$h_i'=\m{h_{N_0-1}}^*\cdots \m{h_0}^*d_i\m{h_0}\cdots
\m{h_{N_0-1}}=\m{h_{N_0-1}}^*\cdots
\m{h_{n_i}}^*\m{uvh_{n_i}}\m{h_{n_i+1}}\cdots \m{h_{N_0-1}}\in C_0(V).$$
For $i=0$, since $d_0h_0\in C_0(V)$, again by \ref{pteo1}
$h_0'=\m{h_{N_0-1}}^*\cdots \m{h_0}^*d_0\m{h_0}\cdots \m{h_{N_0-1}}\in C_0(V)$.
This shows the claim.\fimaf

Define $f'=\m{h_{N_0-1}}^*\cdots\m{h_0}^*f\m{h_0}\cdots\m{h_{N_0-1}}$. By the
fact that $\sigma_{|_{\text{supp}(h_j)}}$ is a homeomorphism it follows that
$f(x_{N_0})=f(x)$. Moreover, since $x_{N_0}\notin V$, by the fact that $V$ is
$\sigma^{-1}$-invariant and $x\notin V$, it follows that $h_i'(x_{N_0})=0$ for
each $i$. Since $f',h_i'\in C(X)$ we have that
$$\|f'-(h_0'+\sum\limits_{i=1}^nh_i')\|_\infty=\|\m{h_{N_0-1}}^*\cdots\m{h_0}^*(
f - ( d _ 0 + \sum\limits_{i=1}^Nd_i))\m{h_0}\cdots\m{h_{N_0-1}}\|\leq\|f-(d_0+\sum\limits_{i=1}^Nd_i)\|<\varepsilon,$$ from where
$|f(x)|=|(f'-(h_0'+\sum\limits_{i=1}^nh_i'))(x_{N_0})|< \varepsilon$.

It remains to analyze the case $x\notin V$, $x\in U$
but $\sigma^n(x)\notin U$ for some $n\leq N_0-1$.
For $i\in\{0,\cdots.n-2\}$ define $h_j$ as above, that is, $h_j\in C_c(U)$
such that $h_j(x_1)=1$, $0\leq h_j \leq 1$ and $\sigma_{|_{\text{supp}(h_j)}}$
is a homeomorphism. For $x_{n-1}$ choose $h_{n-1}\in C_c(U)$ such that $0\leq
h_{n-1}\leq 1$, $h_{n-1}(x_{n-1})=1$, $\sigma_{|_{\text{supp}(h_{n-1})}}$
is a homeomorphism and $\sigma(\text{supp}(h_{n-1}))\subseteq X\setminus K$.
It is possible to choose such $h_{n-1}$ because
$\sigma(x_{n-1})=\sigma^n(x)\in X\setminus U\subseteq X\setminus K$.

\af{ 3}{For $n_i\geq n+1$, $\m{h_{n-1}}^*\cdots
\m{h_0}^*d_i\m{h_0}\cdots \m{h_{n-1}}=0$.}

Denote by $u$ the element $\m{h_{n-2}}^*\cdots
\m{h_0}^*\m{f_1^i}\cdots \m{f_{n-1}^i}$ which is an element of $C(X)$. Then
$$\m{h_{n-1}}^*\cdots\m{h_0}^*\m{f_1^i}\cdots
\m{f_{n+1}^i}=\m{h_{n-1}}^*\m{uf_n^i}\m{f_{n+1}^i}=\me{L(h_{n-1}^*uf_n^i)f_{n+1}^i}.$$
We will show that $L(h_{n-1}^*uf_n^i)f_{n+1}^i=0$.
If $x\notin\supp(f_{n+1}^i)$ or if $\sigma^{-1}(x)=\emptyset$ then \lb$(L(h_{n-1}^*uf_n^i)f_{n+1}^i)(x)=0$.
Suppose therefore $x\in \supp(f_{n+1})$ and $y\in
\sigma^{-1}(x)$. Supposing that $y\in \sigma^{-1}(x)\cap \supp(h_{n-1})$
we have that $x=\sigma(y)\in \sigma(\supp(h_{n-1}))\subseteq
X\setminus K$, which is an absurd because $x\in K$. Therefore if $y\in
\sigma^{-1}(x)$ then $y\notin \supp(h_{n-1})$, and by this way
$L(h_{n-1}^*uf_n^i)(x)=\sum\limits_{y\in \sigma^{-1}(x)}(h_{n-1}^*uf_n^i)(y)=0.$
So $L(h_{n-1}^*uf_n^i)f_{n+1}^i=0$ and the claim is proved.

\af{ 4}{For $n_i\leq n$, $h_i'=\m{h_{n-1}}^*\cdots
\m{h_0}^*d_i\m{h_0}\cdots \m{h_{n-1}}\in C_0(V)$.}

The proof of this claim is analogous to the proof of claim 2.\fimaf

Again $\m{h_{n-1}}^*\cdots \m{h_0}^*f\m{h_0}\cdots
\m{h_{n-1}}=f'$ com $f'(x_n)=f(x)$. Moreover, by the fact that $x_n\notin V$ it follows that
$h_i'(x_n)=0$ for each $i$.
Then
$$\|f'-h_0'-\sum\limits_{n_i\leq
n}h_i'\|=\|\m{h_{n-1}}^*\cdots\m{h_0}^*(f-(d_0+\sum\limits_{i=1}^Nd_i))\m{h_0}
\cdots\m{h_{n-1}}\|< \varepsilon$$ from where
$|f(x)|=|(f'-h_0'-\sum\limits_{n_i\leq n}h_i')(x_n)|< \varepsilon.$

In this way, given $\varepsilon > 0$, for all $x\notin V$, we have that
$|f(x)|\leq\varepsilon$. Therefore $f\in C_0(V)$. \fim

The following theorem is the main result of this section.

\begin{teorema}\label{teobij}
There exists a bijection between the gauge-invariant ideals of $\prcr{X}$ and the open
$\sigma,\sigma^{-1}$-invariant subsets of $X$.
\end{teorema}

\demo All what we have to do is to show that the map
$$\Phi:\{\text{ideals of } \prcr{X}\}\rightarrow \{\text{open
}\sigma,\sigma^{-1}\text{-invariant subsets of } X\},$$ given by $\Phi(I)=V$ where $V$ is
the open subset of $X$ such that $I\cap C(X)=C_0(V)$, is bijective. By the previous
proposition $\Phi$ is surjective. It remains to show that $\Phi$ is injective. For this,
given $I\unlhd \prcr{X}$ gauge-invariant, let $V\subseteq X$ the open subset $\sigma,
\sigma^{-1}$-invariant such that $I\cap C(X)=C_0(V)$. We will show that $\langle
C_0(V)\rangle =I$. It is clear that $\langle C_0(V)\rangle \subseteq I$. By
$\ref{teoiso}$ there exists a *-isomorphism
$\Psi:\frac{\prcr{X}}{\langle C_0(V)\rangle} \rightarrow \mathcal{O}(X',\al',L')$ where
$X'=X\setminus V$. Let $\ovl{\ovl{I}}$ the image of $I$ by the quotient map of
$\prcr{X}$ on $\prcr{X}/\langle C_0(V)\rangle$. Since $\ovl{\ovl{I}}$ is gauge-invariant and $\Psi$ is covariant by the gauge actions we have that $\Psi(\ovl{\ovl{I}})$
is gauge-invariant. Supposing $\ovl{\ovl{I}}\neq 0$, and so $\Psi(\ovl{\ovl{I}})\neq 0$,
it follows that $\Psi(\ovl{\ovl{I}})\cap C(X')=C_0(V')\neq 0$ by
\ref{int.de.id.gauge.inv.com.C(X)}. Let $0\neq g\in C_0(V')$. Then
$g=\Psi(\ovl{\ovl{f}})$ for some $f\in C(X)$ and
$g=\Psi(\ovl{\ovl{a}})$ with $a\in I$. Therefore
$\Psi(\ovl{\ovl{f}})=g=\Psi(\ovl{\ovl{a}})$ from where $\ovl{\ovl{f}}=\ovl{\ovl{a}}$. In
this way, $f-a\in \langle C_0(V)\rangle \subseteq I$ and so $f\in I$. It follows that
$f\in I\cap C(X)=C_0(V)$, that is, $g=\Psi(\ovl{\ovl{f}})=0$, which is an absurd.
Therefore $\ovl{\ovl{I}}=0$ and this shows that $I=\langle C_0(V)\rangle $. \fim

Notice that we have showed that every gauge-invariant idel $I$ of $\prcr{X}$ is of the
form $\langle C_0(V)\rangle $ where $V$ is the $\sigma,\sigma^{-1}$-invariant open subset
such that $I\cap C(X)=C_0(V)$. By this theorem we have the following non simplicity criteria of
$\prcr{X}$:

\begin{corolario}
If $U$ is nonempty and $U\cup \sigma(U)$ is not dense in $X$ then $\prcr{X}$
has at least one gauge-invariant nontrivial ideal.
\end{corolario}

\demo Note that $V=X\setminus\ovl{U\cup\sigma(U)}$ is an open $\sigma,\sigma^{-1}$-invariant set. Since $U\cup
\sigma(U)$ is not dense in $X$ it follows that $V$ is nonempty. Then
$\langle C_0(V)\rangle $ is a nonzero gauge-invariant ideal of $\prcr{X}$. By
the previous theorem, supposing $\langle C_0(V)\rangle =\prcr{X}$ we have
that $C_0(V)=C(X)$, which is a contradiction, because $V\neq X$, by the fact
that $U$ is nonempty. \fim

\section{Topologically free transformations}

In this section we prove that under certain hypothesis about $X$, every ideal of
$\prcr{X}$ has nonzero intersection with $C(X)$ and based on this fact we
show a relationship between the ideals of $\prcr{X}$ and the
$\sigma,\sigma^{-1}$-invariant open subsets of $X$. Also we show a simplicity criteria
for the Cuntz-Krieger algebras for infinite matrices.

\subsection{The theorem of intersection of ideals of $\prcr{X}$ with $C(X)$}

Let us begin with the lemma:

\begin{lema}\label{lsupp}
a)For each $f\in C_c(U)$, $\supp(L(f))\subseteq \sigma(\supp(f))$.\nl b) Let $h,f_1,\cdots
,f_n,g_1,\cdots ,g_n$ be elements of $C_c(U)$ such that $\sigma^{n-1}(\supp(h))\subseteq U$. Then \lb
$\supp(\m{f_k}^*\cdots \m{f_1}^*h\m{g_1}\cdots \m{g_k})\subseteq \sigma^k(\supp(h))$ for
each $k\in \{0,\cdots ,n\}$.
\end{lema}

\demo a) The proof of this fact is similar to the proof given in
[\ref{exeldynamicsystem}:8.7], although our context is a little different. Let $x\in X$
with $L(f)(x)\neq 0$. Suppose $x\notin \sigma(\supp(f))$. Choose $g\in C(X)$ such that
$g(x)=1$ and $g_{|_{\sigma(\supp(f))}}=0$. If $y\in \supp(f)$ then $\al(g)(y)=g(\sigma(y))=0$ because 
$\sigma(y)\in \sigma(\supp(f))$. This shows that
$f\al(g)=0$. So we have $$0\neq L(f)(x)=L(f)(x)g(x)=(L(f)g)(x)=L(f\al(g))=0,$$ which is an absurd.
Therefore $x\in \sigma(\supp(f))$.\nl b) By a) we have that 
$\supp(\m{f_1}^*h\m{g_1})=\supp(L(f_1^*hg_1))\subseteq
\sigma(\supp(f_1^*hg_1)),$ and it is clear that $\sigma(\supp(f_1^*hg_1))\subseteq
\sigma(\supp(h))$. Suppose that
$$\supp(\m{f_{k-1}}^*\cdots \m{f_1}^*h\m{g_1}\cdots \m{g_{k-1}})\subseteq
\sigma^{k-1}(\supp(h)) \text{ for } 2\leq k\leq n.$$ Then, by placing $g=\m{f_{k-1}}^*\cdots
\m{f_1}^*h\m{g_1}\cdots \m{g_{k-1}}$, by a) we have that
$$\supp(\m{f_k}^*g\m{g_k})=\supp(L(f_k^*gg_k))\subseteq
\sigma(\supp(f_k^*gg_k)).$$ Since $\supp(f_k^*gg_k)\subseteq \supp(g)$, and by the
induction hypothesis $\supp(g)\subseteq \sigma^{k-1}(\supp(h))$, it follows that
$\supp(f_k^*gg_k)\subseteq\sigma^{k-1}(\supp(h))$. By hypothesis we have that
$\sigma^{k-1}(\supp(h))\subseteq U$ and so $\sigma(\supp(f_k^*gg_k))\subseteq
\sigma^k(\supp(h))$. \fim

For each $i\neq j$ in $\N$ define $$V^{i,j}=\{x\in X:\,\sigma^i(x)=\sigma^j(x)\}.$$ Note
that for $x\in X$ to be an element of $V^{i,j}$ it is necessary that $x$ lies in
$\text{dom}(\sigma^i)\cap \text{dom}(\sigma^j)$.

\begin{lema}\label{h.que.anula}
If $f_1,\cdots f_i,g_1,\cdots ,g_j\in C_c(U)$ with $i\neq j$ then for each $x\notin
V^{i,j}$ there exists $h\in C(X)$ such that $0\leq h\leq 1$, $h(x)=1$, and
$h\m{f_1}\cdots \m{f_i}\m{g_j}^*..\m{g_1}^*h=0$.
\end{lema}

\demo By taking adjoints we may suppose that $i> j$, and so $i> 0$. Define the set\lb
$K=(\bigcup\limits_{r=1}^i\supp(f_r))(\bigcup\limits_{s=1}^j \supp(g_s))$ which is a
compact subset of $U$. If $x\notin U$, take $h\in C(X)$, $0\leq h\leq 1$ , $h(x)=1$ and
$h_{|_K}=0$. Then $hf_1=0$, which proves the lemma in this case. So we may suppose that
$x\in U$. We will consider two cases: the first when $x\notin \text{dom}(\sigma^i)$ and
the second when $x\in \text{dom}(\sigma^i)$. Suppose $x\notin \text{dom}(\sigma^i)$. Then
there exists $1\leq k\leq i-1$ such that $\sigma^k(x)\notin U$ (note that $i\geq 2$
because $x\in U=\text{dom}(\sigma)$). So $\sigma^k(x)\notin K$. Take $V_0\subseteq X$ an
open subset with $\sigma^k(x)\in V_0$ and $V_0\cap K=\emptyset$. Then
$V=\sigma^{-k}(V_0)\ni x$ is an open subset in $U$. Choose $h\in C_c(U)$ with
$\supp(h)\subseteq V$, $0\leq h\leq 1$ and $h(x)=1$. Then, since
$\sigma^{k-1}(\supp(h^2))\subseteq \sigma^{k-1}(V)\subseteq U$, by $\ref{lsupp}$ b),
$$\supp(\m{f_k}^*\cdots \m{f_1}^*h^2\m{f_1}\cdots \m{f_k})\subseteq
\sigma^k(\supp(h^2))\subseteq \sigma^k(V)\subseteq V_0.$$ Since $V_0\cap K=\emptyset$ and
$\supp(f_{k+1})\subseteq K$ we have that $(\m{f_k}^*\cdots \m{f_1}^*h^2\m{f_1}\cdots
\m{f_k})\m{f_{k+1}}=0$ from where $h\m{f_1}\cdots\m{f_{k+1}}\cdots\m{f_i}=0$. Therefore
$h\m{f_1}\cdots\m{f_i}\m{g_j}^*\cdots\m{g_1}^*h=0$. It remains to show the case $x\in
\text{dom}(\sigma^i)$. By the fact that $i> j$ it follows that $x\in \text{dom}(\sigma^j)$.
Therefore, since $x\notin V^{i,j}$ we have that $\sigma^i(x)\neq \sigma^j(x)$. Let $V_i\ni
\sigma^i(x)$ and $V_j\ni \sigma^j(x)$ open subsets such that $V_i\cap V_j=\emptyset.$ Let
$V=\sigma^{-i}(V_i)\cap \sigma^{-j}(V_j)$ and note that $V$ is an open subset which
contains $x$. Take $h\in C_c(U)$ with $0\leq h\leq 1$, $h(x)=1$ and $\supp(h)\subseteq
V$. Then, since $\sigma^{i-1}(V)\subseteq U$ and $\sigma^{j-1}(V)\subseteq U$, by
$\ref{lsupp}$ b) we have that
$$\supp(\m{f_i}^*\cdots \m{f_1}^*h^2\m{f_1}\cdots \m{f_i})\subseteq
\sigma^i(\supp(h^2))\subseteq V_i$$ and $$\supp(\m{g_j}^*\cdots \m{g_1}^*h^2\m{g_1}\cdots
\m{g_j})\subseteq \sigma^j(\supp(h^2))\subseteq V_j.$$Since $V_i$ and $V_j$ are disjoints
it follows that
$(\m{f_i}^*\cdots \m{f_1}^*h^2\m{f_1}\cdots \m{f_i})(\m{g_j}^*\cdots
\m{g_1}^*h^2\m{g_1}\cdots \m{g_j})=0,$ from where
$h\m{f_1}\cdots \m{f_i}\m{g_j}^*\cdots \m{g_1}^*h=0.$ \fim

\begin{definicao}
We say that the pair $(X,\sigma)$ is topologically free if for each $V^{i,j}$, the
closure $\ovl{V^{i,j}}$ in $X$ has empty interior.
\end{definicao}

By the Baire's theorem, $X$ is topologically free if $\bigcup\limits_{i,j\in
\N}\ovl{V^{i,j}}$ has empty interior. In this way, $Y=X\setminus \bigcup\limits_{i,j\in
\N}\ovl{V^{i,j}}$ is dense in $X$.

Let $S$ be the set of positive linear functionals of $\prcr{X}$ given by
$$S=\{\varphi:\,\varphi \text{ is a positive linear functional and } \varphi_{|_{C(X)}}=\delta_y \text{ for some }y\in Y\}$$
where $\delta_y(f)=f(y)$ for each $f\in C(X)$. We don't know the characteristic of these
functionals, nevertheless for $a\in \prcr{X}$ and $f\in C(X)$ it holds the following
relation:

\begin{lema}\label{flp}
If $\varphi$ is a positive linear functional of $\prcr{X}$ such that
$\varphi_{|_{C(X)}}=\delta_x$ for some $x\in X$ then for each $f\in C(X)$ and $a\in
\prcr{X}$ we have that $\varphi(fa)=\varphi(f)\varphi(a)$ and
$\varphi(af)=\varphi(a)\varphi(f)$.
\end{lema}

\demo By taking adjoints it suffices to prove the case $\varphi(af)=\varphi(a)\varphi(f)$.
For each \lb$b\in \prcr{X}$ we have that
$(b-\varphi(b))^*(b-\varphi(b))\geq 0.$
Therefore if $\varphi$ is a positive functional then 
$\varphi(b^*b)-\varphi(b^*)\varphi(b)=\varphi((b-\varphi(b))^*(b-\varphi(b)))\geq
0,$ from where $\varphi(b)^*\varphi(b)\leq \varphi(b^*b)$. Since $f^*a^*af\leq
f^*f\|a\|^2$ it follows that $\varphi(f^*a^*af)\leq \varphi(f^*f)\|a\|^2$. Put $b=af$,
and so\lb
$0\leq \varphi(af)^*\varphi(af)\leq\varphi(f^*a^*af)\leq \varphi(f^*f)\|a\|^2=\|a\|^2|f(x)|^2,$ where $x$ is
such that $\varphi_{|_{C(X)}}=\delta_x$.
This shows that if $f(x)=0$ then $\varphi(af)=0$. Define $g=f-f(x)$. Then $g(x)=0$ and so
$\varphi(ag)=0$. By this way
$$\varphi(af)-\varphi(a)\varphi(f)=\varphi(af)-\varphi(a)f(x)=\varphi(af)-\varphi(af(x))=\varphi(a(f-f(x)))=\varphi(ag)=0$$
and the lemma is proved.\fim

For each $a\in \prcr{X}$ define
$$\tnorm a\tnorm =\text{sup}\{|\varphi(a)|:\varphi\in S\}$$ which is a seminorm for $\prcr{X}$.

We are not able to show that $\tnorm \text{ }\tnorm $ is nondegenerated in $\prcr{X}$,
but in $L_n$ $\tnorm \text{ }\tnorm $ has the property, given by the following lemma,
that $\tnorm r\tnorm \neq 0$ for every positive nonzero element of $L_n$,
remembering that $L_n=C(X)+K_1+\cdots +K_n$ for each $n\geq 1$ and $L_0=C(X)$.

\begin{lema}\label{ndeg}
Let $(X,\sigma)$ be topologically free. For each $r\in L_n$ with $r\geq 0$ and $r\neq 0$
it holds that $\tnorm r\tnorm \neq 0$.
\end{lema}

\demo

\vspace{0cm}\af{ 1}{If $0\neq r\in L_n$, $r$ positive and $r\notin C(X)$ then there
exists $g\in C_c(U)$ with $\sigma_{|_{\supp(g)}}$ a homeomorphism and $\m{g}^*r\m{g}\neq
0$}

Since $r\geq 0$ we may write $r=b^*b$ with $b\in L_n$. Suppose that for each $g\in
C_c(U)$ with $\sigma_{|_\supp(g)}$ homeomorphism, it holds that $\m{g}^*r\m{g}=0$, and so
$\m{g}^*b^*=0$. Then (making use of partition of unity we may write each $f\in C_c(U)$ as a
sum of $g$ as above) we have that $\m{f}^*b^*=0$ for each $f\in C_c(U)$ and so
$M^*b^*=0$. It follows that $K_1b^*=0$, and since $C_0(U)\subseteq K_1$ by $\ref{l1}$ b) we
have that $C_0(U)b^*=0$ and by $\ref{multired}$ b) it follows that $b^*\in C(X)$. In this
way $r=b^*b\in C(X)$, which contradicts the hypothesis and the claim is proved.

\af{ 2}{If $0\neq r\in L_n$, $r\geq 0$ and $r\notin C(X)$ then there exists
$g_1,\cdots,g_i\in C_c(U)$ such that $\sigma_{|_{\supp(g_j)}}$ is a homeomorphism for
each $j$ and
$0\neq \m{g_i}^*\cdots \m{g_1}^*r\m{g_1}\cdots \m{g_i}\in
C(X).$}

By claim 1 there exists $g_1\in C_c(U)$ such that $\sigma_{|_{\text{supp}(g_1)}}$ is
homeomorphism and $0\neq \m{g_1}^*r\m{g_1}$. Note that $\m{g_1}^*r\m{g_1}\in L_{n-1}$. By induction suppose $0\neq
\m{g_l}^*\cdots\m{g_1}^*r\m{g_1}\cdots\m{g_l}\in L_1$ where $g_j\in C_c(U)$ and
$\sigma_{|_{\text{supp}(g_j)}}$ is a homeomorphism for each $j$. Then, by claim 1, or
$\m{g_l}^*\cdots\m{g_1}^*r\m{g_1}\cdots\m{g_l}\in C(X)$ or there exists $g_{l+1}\in
C_c(U)$ with $\sigma_{|_{\text{supp}(g_{l+1})}}$ homeomorphims and $0\neq
\m{g_{l+1}}^*\m{g_l}^*\cdots\m{g_1}^*r\m{g_1}\cdots\m{g_l}\m{g_{l+1}}.$ Since
$\m{g_{l+1}}\m{g_l}^*\cdots\m{g_1}^*r\m{g_1}\cdots\m{g_l}\m{g_{l+1}}\in C(X)$
the claim is proved.\fimaf

We will now show the lemma. Let $r\in L_n$, $r$ positive and no
null. It is enough to show that there exists $\varphi\in S$ such that $\varphi(r)\neq 0$.
Since $(X,\sigma)$ is topologically free then $Y\,\,(=X\setminus
\bigcup\limits_{i,j}\ovl{V^{i,j}})$ is dense in $X$. So, if $r\in C(X)$ then there exists
$y\in Y$ such that $r(y) > 0$. Take $\varphi$ which extends $\delta_y$, and therefore
$\varphi(r)\neq 0$. Suppose $r\notin C(X)$. Choose $f_{x_1},\cdots ,f_{x_i}\in C_c(U)$
as in claim 2. Then $0\neq h=\m{f_{x_i}}^*\cdots \m{f_{x_1}}^*r\m{f_{x_1}}\cdots
\m{f_{x_i}}\in C(X).$ So
$$h^*hh^*=\m{f_{x_i}}^*\cdots \m{f_{x_1}}^*r\m{f_{x_1}}\cdots
\m{f_{x_i}}h\m{f_{x_i}}^*\cdots \m{f_{x_1}}^*r\m{f_{x_1}}\cdots \m{f_{x_i}}\neq 0$$ from
where
$g=\m{f_{x_1}}\cdots \m{f_{x_i}}h\m{f_{x_i}}^*\cdots \m{f_{x_1}}^*\neq 0.$ How
$\sigma_{|_{\supp(f_{x_i})}}$ is homeomorphism it follows by $\ref{l1}$ a) that
$\m{f_{x_i}}h\m{f_{x_i}}^*\in C(X).$ Applying these arguments sucessively it may be proved
that \lb$g=\m{f_{x_1}}\cdots \m{f_{x_i}}h\m{f_{x_i}}^*\cdots \m{f_{x_1}}^*\in C(X).$ By
the the same argments it follows that\lb
$u=\m{f_{x_1}}\cdots \m{f_{x_i}}\m{f_{x_i}}^*\cdots \m{f_{x_1}}^*\in C(X).$ Since $g\neq
0$ there exists $y\in Y$ such that $g(y)\neq 0$. Take $\varphi\in S$ which extends
$\delta_y$. Then we have that $\varphi(g)=g(y)\neq 0$. By $\ref{flp}$, since $g=uru$,
$\varphi(g)=\varphi(uru)=\varphi(u)\varphi(r)\varphi(u)$ and therefore
$\varphi(r)\neq 0.$ \fim

Now we are able to prove the main result of this section.

\begin{teorema}\label{top livre}
If $(X,\sigma)$ is topologically free then each nonzero ideal of
$\prcr{X}$ has nonzero intersection with $C(X)$.
\end{teorema}

\demo By $\ref{int ideal com C(X)}$ it suffices to prove that every nonzero
ideal of $\prcr{X}$ has nonzero intersection with $K$. Let $0\neq I\unlhd
\prcr{X}$. Suppose $I\cap K=0$. Then the quotient *-homomorphism
$\pi:\prcr{X}\rightarrow \prcr{X}/I$ is such that $\pi_{|_K}$
is an isometry.

\af{}{For each $b\in \prcr{X}$ it holds that
$\tnorm E(b)\tnorm \leq \|\pi(b)\|$ where $E$ is the conditional expectation defined in section \ref{ck} }.

Let $a$ be of the form $a=\sum\limits_{0\leq i\leq n\atop 0\leq j\leq m} a_{i,j}$ with
$a_{0,0}\in C(X)$ and $a_{i,j}\in M^i{M^j}^*$ for $i\neq 0$ or $j\neq 0$,
$a_{i,j}=\sum\limits_{1\leq k\leq n_{i,j}} a_{i,j}^k$, $a_{i,j}^k=\m{f_{i,j,1}^k}\cdots
\m{f_{i,j,i}^k}\m{g_{i,j,1}^k}^*\cdots \m{g_{i,j,j}^k}^*$ where
$f_{i,j,l}^k,g_{i,j,t}^k\in C_c(U)$ for each $i,j,k,l$ and $t$. Given $\varepsilon > 0$
there exists $\varphi\in S$ which extends $\delta_y$ for some $y\in Y$ such that
$\tnorm E(a)\tnorm -\varepsilon \leq |\varphi(E(a))|.$ Note that $y\notin V^{i,j}$ for $i\neq
j$. Then, for every $a_{i,j}^k$ with \lb$i\neq j$, by \ref{h.que.anula} there exists
$h_{i,j}^k\in C(X)$, $0\leq h_{i,j}^k\leq 1$, such that $h_{i,j}^k(y)=1$ and
$ha_{i,j}^kh=0$. Define
$h=\prod\limits_{0\leq i\leq n\atop 0\leq j\leq m}\prod\limits_{1\leq k\leq
n_{i,j}}h_{i,j}^k.$ Then $ha_{i,j}h=0$ for each $i\neq j$ from where $hah=hE(a)h$, and
moreover $h(y)=1$. By $\ref{flp}$
$\varphi(hE(a)h)=\varphi(h)\varphi(E(a))\varphi(h)=h(y)\varphi(E(a))h(y)=
\varphi(E(a)),$ and so $\varphi(E(a))=\varphi(hE(a)h)=\varphi(hah).$ Since
$hah=hE(a)h\in K$ e $\pi_{|_K}$ is an isometry it follows that $\|hah\|=\|\pi(aha)\|$.
Then
$$\tnorm E(a)\tnorm -\varepsilon\leq |\varphi(E(a))|=|\varphi(hah)|\leq
\|hah\|=\|\pi(hah)\|\leq\|\pi(a)\|.$$ Since $\varepsilon$ is arbitrary it follows that
$\tnorm E(a)\tnorm \leq \|\pi(a)\|$ for $a$ in this form. Given $b\in \prcr{X}$, for each
$\varepsilon > 0$ choose $a\in \prcr{X}$ as above such that $\|a-b\|\leq \varepsilon$.
Then
$$\tnorm E(b)\tnorm \leq \tnorm E(b-a)\tnorm +\tnorm E(a)\tnorm \leq \tnorm E(a)\tnorm +\varepsilon\leq
\|\pi(a)\|+\varepsilon\leq$$ $$\leq\|\pi(a-b)\|+\|\pi(b)\|+\varepsilon\leq \|\pi(b)\|+2\varepsilon.$$
Again, since $\varepsilon$ is arbitrary it follows that $\tnorm E(b)\tnorm \leq\|\pi(b)\|$,
and the claim is proved.\fimaf

Observe that $\ovl{E(I)}$ is a closed ideal of $K$. Also, $\ovl{E(I)}$ is
nonzero, because $0\neq I$ and $E$ is faithful. Then $\ovl{E(I)}\cap L_n\neq 0$ for
some $n$ (see [\ref{davidson}: III.4.1]). Let $0\neq c\in \ovl{E(I)}\cap L_n$.
Then, since $c^*c\in L_n$ and $c^*c$ is positive and nonzero it follows by
$\ref{ndeg}$ that $\tnorm c^*c\tnorm \neq 0$. We shall prove that $\tnorm c^
*c\tnorm =0$, and this will be an absurd. For each $a=E(b)\in E(I)$ with $b\in
I$ we have that $$\tnorm a^*a\tnorm =\tnorm E(b^*)E(b)\tnorm =\tnorm
E(b^*E(b))\tnorm \leq \|\pi(b^*E(b))\|.$$ By the fact that $b^*E(b)\in I$ it follows that
$\pi(b^*(E(b)))=0$ and so $\tnorm a^*a\tnorm =0$. This shows that $\tnorm
a^*a\tnorm =0$ for each $a\in E(I)$. Given $\varepsilon > 0$, take $a\in E(I)$
such that $\|a^*a-c^*c\|\leq \varepsilon$. Then $$\tnorm c^*c\tnorm \leq \tnorm
c^*c-a^*a\tnorm +\tnorm a^*a\tnorm =\tnorm c^*c-a^*a\tnorm \leq
\|c^*c-a^*a\|\leq\varepsilon.$$ So $\tnorm c^*c\tnorm \leq\varepsilon$ for each
$\varepsilon > 0$ from where $\tnorm c^*c\tnorm =0$, and that is an absurd.
Therefore $I\cap K\neq 0$, and the theorem is proved. \fim

\subsection{Relationship between the ideals of $\prcr{X}$ and the
$\sigma,\sigma^{-1}$-invariant open subsets of $X$}

We obtain here a relationship between the ideals of $\prcr{X}$ and the
$\sigma,\sigma^{-1}$-invariant open subsets of $X$ under an additional hypothesis about
$(X,\sigma)$, which is that for every closed $\sigma,\sigma^{-1}$-invariant subset $X'$ of $X$, $(X',\sigma_{|_{X'}})$ is topologically free.

\begin{proposicao}\label{id.C0(V)}
Let $I$ be an ideal of $\prcr{X}$ and $V\subseteq X$ the open subset such that \lb$I\cap
C(X)=C_0(V)$. If $(X',\sigma_{|_{X'}})$ is topologically free (where $X'=X\setminus V$)
then $I=\langle C_0(V)\rangle $. \end{proposicao}

\demo By \ref{sigma invariante} $V$ is $\sigma,\sigma^{-1}$-invariant, from where
$X'$ is also $\sigma,\sigma^{-1}$-invariant. By \ref{teoiso} there exists a *-isomorphism
$\Psi:\frac{\prcr{X}}{\langle C_0(V)\rangle }\rightarrow \mathcal{O}(X',\al',L')$. Obviously $\langle C_0(V)\rangle \subseteq I$. Suppose
$I\neq \langle C_0(V)\rangle $. Then $\ovl{\ovl{I}}\neq 0$ and so
$\Psi(\ovl{\ovl{I}})\neq 0$. By \ref{top livre}, $\Psi(\ovl{\ovl{I}})\cap
C(X')\neq 0$. Let $0\neq g\in \Psi(\ovl{\ovl{I}})\cap C(X')$. Then $g=\Psi(\ovl{\ovl{a}})$ for some $a\in I$
and also $g=\Psi(\ovl{\ovl{f}})$, because $\Psi(\ovl{\ovl{C(X)}})=C(X')$. Therefore
$\Psi(\ovl{\ovl{a}})=\Psi(\ovl{\ovl{f}})$ from where $\ovl{\ovl{a}}=\ovl{\ovl{f}}$ and so
$f-a\in \langle C_0(V)\rangle \subseteq I$, in other words, $f\in I$. In this way $f\in
I\cap C(X)=C_0(V)$ and so $\ovl{\ovl{f}}=0$ from where $g=\Psi(\ovl{\ovl{f}})=0$, which is
a absurd. So we conclude that $I=\langle C_0(V)\rangle $. \fim

\begin{teorema}\label{teoisoideais}
If $(X,\sigma)$ is such that $(X',\sigma_{|_{X'}})$ is topologically free for every
closed subset \lb$\sigma,\sigma^{-1}$-invariant $X'$ of $X$ then every ideal of $\prcr{X}$ is
of the form $\langle C_0(V)\rangle $ for some open subset $V\subseteq X$. Moreover, the map
$V\longrightarrow \langle C_0(V)\rangle$ is a bijection 
between the open $\sigma,\sigma^{-1}$-invariante subsets of $X$ and the ideals of $\prcr{X}$.
\end{teorema}

\demo Let $I\unlhd \prcr{X}$, and $C_0(V)=I\cap C(X)$. By \ref{sigma invariante} $V$
is $\sigma,\sigma^{-1}$-invariant, from where $X'=X\setminus V$ is also
$\sigma,\sigma^{-1}$-invariant. By hypothesis $(X',\sigma_{|_{X'}})$ is topologically
free. By \ref{id.C0(V)}, $I=\langle C_0(V)\rangle $. In particular, note that every ideal of $\prcr{X}$ is
gauge-invariant. So, by \ref{teobij} the map a $V\longrightarrow \langle C_0(C)\rangle$ is a bijection. \fim

\subsection{An simplicity criteria for the Cuntz-Krieger algebras for infinite matrices}

Recall that $G_R(A)$ is the oriented graph whose vertex are the elements of $G$ such
that given $x,y\in G$ there exists an oriented edge from $x$ to $y$ if $A(x,y)=1$. An
path from $x$ to $y$ is a finite sequence $x_1\cdots x_n$ such that $x_1=x$, $x_n=y$ and
$A(x_i,x_{i+1})=1$ for each $i$. We will say that $G_R(A)$ é transitive if for each
$x,y\in G$ there exists a path from $x$ to $y$.

The following proposition singles out the $\sigma,\sigma^{-1}$-invariant open subsets of
$\m{\Omega_A}$.

\begin{proposicao}\label{carac.abertos.sigma.inv.de.OA}
If $G_R(A)$ is transitive, the unique $\sigma$-invariants nonempty open subsets of
$\widetilde{\Omega_A}$ are $\widetilde{\Omega_A}\setminus \emptyset$ and
$\widetilde{\Omega_A}$.
\end{proposicao}

\demo Let $V$ be a $\sigma$-invariant open subset of $\widetilde{\Omega_A}$. Let $\xi\in
V$ an element whose stem is infinite. (such elements form a dense subset in
$\widetilde{\Omega_A}$). Choose $V_n$ neighbourhood of $\xi$ in $V$,
$$V_n=\{\nu\in \widetilde{\Omega_A} ;
w(\nu)_{|_n}=w(\xi)_{|_n}\}$$ where $w(\nu)$ is the stem of $\nu$. Let
$\mu\in\widetilde{\Omega_A}$ such that $|w(\mu)|\geq 1$ and let $x\in G$, with $x\in
\mu$. Since $G_R(A)$ is transitive there exists a path $x_1\cdots x_m$ from $w(\xi)_n$ to
$x$, and by this way\lb $w(\xi)_{|_n}x_2\cdots x_{m-1}\mu\in V_n\subseteq V$. Since $V$ is
$\sigma$-invariant it follows that $\mu\in V$ because \lb$\mu=\sigma^{n+m-2}(w(\xi)_{|_n}x_2\cdots
x_{m-1}\mu)$. So $U\subseteq V$. If $\emptyset \neq \xi\in \m{\Omega_A}\setminus
U$ then there exists $x\in G$ such that $x^{-1}\in \xi$. Since $x\xi\in U\subseteq V$ and
$\sigma(x\xi)=\xi$ it follows that $\xi\in V$. This shows that $\m{\Omega_A}\setminus
\emptyset\subseteq V$, from where the result follows. \fim

Since $\m{\Omega_A}$ and $\m{\Omega_A}\setminus \emptyset$ are $\sigma^{-1}$-invariant it
follows by the previous proposition that the unique $\sigma,\sigma^{-1}$-invariant open
nonempty subsets of $\m{\Omega_A}$ are $\m{\Omega_A}$ and $\m{\Omega_A}\setminus
\emptyset$.

Given $\xi\in \text{dom}(\sigma^i)$ with $w(\xi)=x_1x_2\cdots $ we have that
$w(\sigma^i(\xi))=x_{i+1}x_{i+2}\cdots $. This shows that if $\xi\in V^{i,j}$ then
$w(\xi)$ is infinite, because if we suppose that $|w(\xi)|=n$, then we have that
$n-i=|w(\sigma^i(\xi))|=|w(\sigma^j(\xi))|=n-j$ from where $i=j$, which is an absurd.

The following proposition shows a relationship between $Gr(A)$ and $\m{\Omega_A}$.

\begin{proposicao}\label{transitivo implica top. livre}
If $Gr(A)$ is transitive then $\m{\Omega_A}$ is topologically free.
\end{proposicao}

\demo Suppose $i > j$, $i=j+k$ and that $\ovl{V^{i,j}}$ has nonempty interior. Let $\nu$
be an interior point of $\ovl{V^{i,j}}$ and $V_\nu\subseteq \ovl{V^{i,j}}$ an open subset
which contains $\nu$. Then there exists an element $\xi\in V_\nu\cap V^{i,j}$. Since
$\sigma^i(\xi)=\sigma^j(\xi)$ we have that
$$x_{i+1}x_{i+2}\cdots =w(\sigma^i(\xi))=w(\sigma^j(\xi))=x_{j+1}x_{j+2}\cdots ,$$ from
where $x_{i+r}=x_{j+r}$ for $r\geq 1$. Since $i=j+k$ it follows that
$x_{i+k}=x_{j+k}=x_i$, and also that
$x_{i+(k+r)}=x_{j+(k+r)}=x_{(j+k)+r}=x_{i+r}$ for each $r\geq 1$. Applying the last
equality repeatedly it follows that $x_{i+nk+r}=x_{i+r}$ for each $n\in \N$ and $r\geq
1$. This shows that \lb$w(\xi)=x_1\cdots x_{i-1}sss\cdots,$ where $s=x_ix_{i+1}\cdots
x_{i+(k-1)}$. Since $w(\xi)$ is infinite, there exists $n\geq i$ such that $V_n=\{\eta\in \m{\Omega_A}:w(\eta)_{|_n}=x_1\cdots x_n=w(\xi)_{|_n}\}\subseteq
V_\nu$.

\af{}{$V_n=\{\xi\}$}

Supposing $\eta\in V_n\cap V^{i,j}$, with the same arguments as above it may be proved
that\lb
$w(\eta)=x_1x_2\cdots x_{i-1}sss\cdots,$
from where $w(\eta)=w(\xi)$, and since $\eta,\xi$ have infinite stems it follows that
$\eta=\xi$. Let $\nu\in V_n$. Then, since $V_n\subseteq \ovl{V^{i,j}}$ there exists a net
$(\nu_l)_l\subseteq V^{i,j}$ such that $\nu_l\rightarrow \nu$. Since $\nu\in V_n$ and
$V_n$ is open we may suppose that $(\nu_l)_l\subseteq V_n$. Therefore $\nu_l=\xi$ for
each $l$ and so $\nu=\xi$. This proves the claim. \fimaf

Let $y\in G\setminus\{x_i,x_{i+1},\cdots ,x_{i+(k-1)}\}$. By the fact that $Gr(A)$ is transitive there
exists a path $y_1\cdots y_r$ where $y_1=x_{n+1}$ and $y_r=y$ and an other path
$z_1\cdots z_t$ such that $z_1=y$ e $z_t=x_1$. In this way we may consider the infinite
admissible word
$x_1\cdots x_ny_1\cdots y_rz_2\cdots z_{t-1}w(\xi)$ which is the stem of some element
$\mu\in\m{\Omega_A}$. Notice that $\mu\in V_n$ by the definition of $V_n$ and that
$\mu\neq \xi$, because its stems are distinct. This contradicts the claim. Therefore,
$\ovl{V^{i,j}}$ has empty interior, and so $\m{\Omega_A}$ is topologically free. \fim

We will prove now the main result of this section.

\begin{proposicao}
If $Gr(A)$ is transitive the unique ideals of $\m{O_A}$ are the null ideal, $O_A$ and
$\m{O_A}$.
\end{proposicao}

\demo By \ref{carac.abertos.sigma.inv.de.OA} the unique closed
$\sigma,\sigma^{-1}$-invariants subsets of $\m{\Omega_A}$ are $\m{\Omega_A}$, the set
$\{\emptyset\}$ (if $\emptyset\in \m{\Omega_A}$, that is, if $O_A\neq \m{O_A}$ by
[\ref{exelmatinf}: 8.5]) and the empty set. Since these subsets are topologically free,
by \ref{teoisoideais} the ideals of $\prcr{\m{\Omega_A}}$ are precisely $0$, $\langle
C_0(\m{\Omega_A}\setminus \emptyset)\rangle $ and $\prcr{\m{\Omega_A}}$. Therefore if
$\emptyset\notin \m{\Omega_A}$ (that is, if $O_A =\m{O_A}$) then $\prcr{\m{\Omega_A}}$
has no nontrivial ideals and the proposition is proved in this case. If $\emptyset\in
\m{\Omega_A}$ then by \ref{teoisoideais} $\prcr{\m{\Omega_A}}$ has exactly one
nontrivial ideal, which is $\langle C_0(\m{\Omega_A}\setminus \emptyset)\rangle $.
Therefore $\m{O_A}$ has also exactly one nontrivial ideal. By [\ref{exelmatinf}: 8.5]
$O_A\neq \m{O_A}$ and since $0\neq O_A\unlhd \m{O_A}$ it follows that $O_A$ is a
nontrivial ideal of $\m{O_A}$, and so is unique. \fim

A direct consequence of this proposition is that if $G_R(A)$ is transitive then $O_A$ is
simple.

\vspace{2cm}
Departamento de Matemática, Universidade Federal de Santa Catarina, Brasil.\nl
E-mail: royer@mtm.ufsc.br\nl\nl 
Departamento de Matemática, Universidade Federal de Santa Catarina, Brasil.\nl
E-mail: exel@mtm.ufsc.br


\addcontentsline{toc}{section}{References}
\begin{thebibliography}{99}

\bibitem{pedersen}\label{pedersen} G. K.
Pedersen, {\it $C^*$-algebras and Their Automorphism Groups}, Academic Press, 1979.

\bibitem{davidson}\label{davidson} K. R. Davidson, {\it $C^*$-Algebras by
Example}, Fields Institute Monographs, 1996.

\bibitem{pimsner}\label{pimsner} M. V. Pimsner, {\it A class of $C^*$-Algebras
generalizing both Cuntz-Krieger Algebras and crossed products by $\Z$}, In Free
probability theory (Waterloo, ON, 1995), volume 12 of Fields Inst. Commun., pages
189-212. Amer. Math. Soc., Providence, RI, 1997.

\bibitem{mclanahan}\label{mclanahan} K. McClanahan, {\it K-theory for partial actions by discrete groups}, J.
Funct. Anal. {\bf 130}, (1995), 77-117.

\bibitem{exelmatinf}\label{exelmatinf} R. Exel and M. Laca, {\it Cuntz-Krieger
Algebras for Infinite Matrices}, J. reine angew. Math. {\bf 521} (1999), 119-172.

\bibitem{exelprodcruz}\label{exelprodcruz} R. Exel, {\it A New Look at The
Crossed-Product of a $C^*$-algebra by an Endomorphism},Ergodic Theory Dynam. Systems, to
appear.

\bibitem{exelprodcruzpartial}\label{exelprodcruzpartial} R. Exel, {\it Circle actions on C$^*$-algebras,
partial automorphisms and a generalized Pimsner-Voiculescu exact sequence}, J. Funct. Anal. {\bf 122}, (1994), 361-401.

\bibitem{exeldynamicsystem}\label{exeldynamicsystem} R. Exel and A. Vershik,
{\it $C^*$-algebras of Irreversible Dynamical Systems},
http://www.arxiv.org/abs/math.OA/0203185

\bibitem{katsura}\label{katsura} T. Katsura, {\it A construction of
$C^*$-algebras from $C^*$-correspondences}, Advances in Quantum Dynamics, 173-182,
Contemp. Math, 335, Amer. Math. Soc., Providence, RI, 2003.

\bibitem{kwasniewski}\label{kwasniewski} B. K. Kwasniewski, {\it Covariance algebra of a partial dynamical
system}, http://lanl.arxiv.org/abs/math.OA/0407352

\end{thebibliography}
\end{document}